\documentclass[11pt]{amsart}
\usepackage{amsmath,amssymb,amsthm}
\usepackage[latin1]{inputenc}
\usepackage{version,tabularx,multicol}
\usepackage{graphicx,float,psfrag}
\usepackage{stmaryrd}

\newcommand{\cT}{{\mathfrak T}}

\newcommand{\reff}[1]{(\ref{#1})}

\theoremstyle{plain}
\newtheorem{theo}{Theorem}[section]

 \newtheorem{corollary}[theo]{Corollary}

 \newtheorem{proposition}[theo]{Proposition}
\newtheorem{lem}[theo]{Lemma}
 \newtheorem{lemma}[theo]{Lemma}
\newtheorem{defi}[theo]{Definition}
 
\theoremstyle{remark}

 \newtheorem{remark}[theo]{Remark}

 \newtheorem{example}[theo]{Example}

\newcommand{\cb}{{\mathcal B}}

\newcommand{\ck}{{\mathcal K}}

\newcommand{\cn}{{\mathcal N}}
\newcommand{\cm}{{\mathcal M}}

\newcommand{\ct}{{\mathcal T}}
\newcommand{\cx}{{\mathcal X}}

\newcommand{\cz}{{\mathcal Z}}

\newcommand{\mE}{{\mathbb E}}

\newcommand{\mN}{{\mathbb N}}

\newcommand{\mP}{{\mathbb P}}
\newcommand{\mQ}{{\mathbb Q}}
\newcommand{\mR}{{\mathbb R}}

\newcommand{\T}{{\cal T}}

\newcommand{\bm}{\mathbf m}
\newcommand{\bN}{{\mathbf N}}
\newcommand{\bP}{{\mathbf P}}

\newcommand{\ind}{{\bf 1}}

\newcommand{\Br}{{\rm Br}}
\newcommand{\Lf}{{\rm Lf}}
\newcommand{\Sk}{{\rm Sk}}

\newcommand{\expp}[1]{\mathop {\mathrm{e}^{ #1}}}

 \def\beqlb{\begin{eqnarray}}\def\eeqlb{\end{eqnarray}}
 \def\beqnn{\begin{eqnarray*}}\def\eeqnn{\end{eqnarray*}}

 \def\ar{\!\!&}

 \def\mbb{\mathbb}

 \def\qed{\hfill$\Box$\medskip}

\newcommand{\bcen}{\begin{center}}
\newcommand{\ecen}{\end{center}}
\newcommand{\bgeqn}{\begin{equation}}
\newcommand{\edeqn}{\end{equation}}

\def\dz{\delta}
\def\lz{\lambda}

\def\ez{\epsilon}

\def\M{{\mathcal  M}}

\def\cn{{\mathcal  N}}

\def\ct{{\mathcal  T}}
\def\Z{{\mathcal  Z}}
\def\mN{{\mbb  N}}
\def\N{{\mbb  N}}
\def\mE{{\mbb  E}}
\def\E{{\mbb  E}}
\def\mP{{\mbb  P}}
\def\mR{{\mbb  R}}

\def\mT{{\mbb  T}}

\def\rar{\rightarrow}

\def\ra{\rangle}
\def\la{\langle}

\def\l{\left}
\def\r{\right}
\newcommand{\lb}{\llbracket}
\newcommand{\rb}{\rrbracket}
\newcommand{\cal}{\mathcal}
 \def\ar{\!\!\!&}

\begin{document}

\title[Pruning and admissible family]{A tree-valued Markov process associated with an admissible
family of branching mechanisms}
\date{\today}

\author{Hongwei Bi}

\address{
Hongwei Bi, School of Insurance and Economics, University of International
Business and Economics, Beijing 100029, P.R.CHINA. }

\email{bihw@uibe.edu.cn}

\author{Hui He}

\address{
Hui He, Laboratory of Mathematics and Complex Systems, School of Mathematical
Sciences, Beijing Normal University, Beijing 100875, P.R.CHINA. }

\email{hehui@bnu.edu.cn}


\begin{abstract}
By studying an admissible family of branching mechanisms introduced in Li
(2014), we obtain a pruning procedure on L\'evy trees. Then we construct
a decreasing L\'evy-CRT-valued process $\{\T_t\}$ by pruning L\'evy trees
and an analogous process $\{\T^*_t\}$ by pruning a critical L\'evy tree
conditioned to be infinite. Under a regular condition on the admissible
family of branching mechanisms, we show that the law of $\{\T_t\}$ at the
ascension time can be represented by $\{\T^*_t\}$. The results generalize
those studied in Abraham and Delmas (2012).

\end{abstract}



\maketitle

\section{Introduction}\label{sec:introduction}
A general pruning procedure was  introduced in Abraham et al.\cite{ADV10}
on L\'evy trees and was further explored by Abraham and Delmas  \cite{AD12}.
In particular,  a decreasing continuum-tree-valued process was constructed
and studied in \cite{AD12} which is associated with a family of branching
mechanisms obtained by shifting a branching mechanism. More precisely,
let $\psi$ be a branching mechanism defined by
\begin{equation}
   \label{eq:psi}
\psi(\lambda)=b\lambda+c\lambda^2+\int_{(0,\infty)}
\left(\expp{-\lambda z}-1+\lambda z\right)m(dz),\quad \lz\geq0,
\end{equation}
where $b\in \mR$, $c\ge 0$ and $m$ is  a $\sigma$-finite measure
on $(0,+\infty)$ such that $\int_0^{\infty}(z\wedge z^2)m(dz)<+\infty$. Define
$\psi^{\theta}(\lz)=\psi(\theta+\lz)-\psi(\theta)$. Denote by $\Theta^{\psi}$
the set of $\theta$ such that $\int_1^{\infty}\expp{-\theta z}m(dz)<\infty$.
The family of branching mechanisms $\{\psi^{\theta}, \theta\in \Theta^{\psi}\}$
was considered in \cite{AD12}.

Li  \cite{Li14} introduced the  {\it admissible family} of branching mechanisms
which generalized those used in \cite{AD12}. Roughly, the model is described as
follows: Given a time interval $\cT\subset\mR$, let $(\theta,\lz)\mapsto\zeta_{\theta}(\lz)$
be a continuous function on $\cT\times[0,\infty)$ with representation
\[
\zeta_{\theta}(\lz)=\beta_{\theta}\lz+\int_{(0,\infty)}(1-\expp{-z\lz})n_{\theta}(dz),\quad
\theta\in\cT,\,\lz\geq0,
\]
where $\beta_{\theta}\geq0$ and $(1\wedge z)n_{\theta}(dz)$ is a finite kernel
from $\cT$ to $(0,\infty)$. Then $\{\psi_{\theta}, \theta\in \cT\}$ is called
an admissible family  if
\[
\psi_{q}(\lz)=\psi_t({\lz})+\int_{t}^{q}\zeta_{\theta}(\lz)d\theta,
\quad q\geq t\in\cT,\, \lz\geq0.
\]
In particular,  $\{\psi^{\theta}, \theta\in \Theta^{\psi}\}$ considered in
\cite{AD12} is an admissible family with
 \[
  \zeta_{\theta}(\lz)=2c\lz+\int_{(0,\infty)}(1-\expp{-z\lz} )z\expp{-z\theta}m(dz).
  \]
  By using the techniques of stochastic equations and measure-valued processes,
  Li \cite{Li14} studied a class of increasing path-valued Markov processes
  associated with the admissible family. Those path-valued processes can be
  regarded as counterparts of the tree-valued processes constructed in \cite{AD12}
  (However, to the best of our knowledge, no link is actually pointed out between
  tree-valued processes and path-valued branching processes).  It is natural to ask
  whether there exists  a continuum-tree-valued process associated with a given
  admissible family by pruning L\'evy trees.

The second motivation of the present work is the study of the so-called
{\it ascension process}. It was first introduced in the  pioneer work of
Aldous and Pitman \cite{AP98}, where they constructed a tree-valued Markov
process $\{ {\cal G}(u)\}$ by pruning Galton-Watson trees (edge percolation
on trees) and an analogous process $\{ {\cal G}^*(u)\}$  by
pruning a critical or subcritical Galton-Watson tree conditioned to be infinite.
It was shown in \cite{AP98} that the process $\{{\cal G}(u)\}$ run until its ascension
time (the first time for which the total mass is finite)
has a representation in terms of $\{{\cal G}^*(u)\}$ in the special case
of Poisson offspring distributions. By using the pruning procedure defined in
 \cite{ADV10} and exploration processes introduced in \cite{LL98a},
 Abraham and Delmas  \cite{AD12} extended the above results to L\'evy trees,
 where a decreasing L\'evy-tree-valued process $\{\T_{\theta}, \theta\in\Theta^{\psi}\}$
 was constructed such that $\T_{\theta}$ is a $\psi^{\theta}$-L\'evy tree.
 They also showed that  $\{\T_{\theta}\}$ run until its ascension time can be
 represented in terms of another tree-valued process obtained by applying the
 same pruning procedure to a L\'evy tree conditioned on non-extinction.
 Similar results can be found in \cite{ADH12a} for Galton-Watson trees where
 the trees are pruned based on bond percolation. The cases for
 sub-trees of L\'evy trees were also studied in \cite{ADH12b}.

In this paper the framework is locally compact measured rooted  real tree
 $(\T, d, \emptyset, \bm)$. The collection is denoted by $\mT$.
Based on the pruning procedure of \cite{AD12, ADV10}, we introduce a more
general pruning mechanism as follows. Let $\T\in\mT, t\in\cT$
and $q\in \cT_t=\cT\cap [t,\infty)$. Put marks on $\cT$ with a
Poisson point measure $M_t^{\T}([t, q], dy)$ as follows:
\begin{enumerate}
\item Assign marks to the skeleton of $\T$ according to a Poisson point
measure with intensity $\int_t^q\beta_{\theta}d\theta \ell^{\T}(dy)$,
where $\ell^{\T}$ is the length measure of $\T$;
\item Assign marks to each node $y\in\T$ of infinite degree with
probability $1-m_{\Delta_y}(t,q)$ (see (\ref{Pmark})).
\end{enumerate}
We prune $\T$ according to the marks and consider the pruned tree $\T_q^t$
containing the root. Theorem 4.2 then gives the connection between
the pruning mechanism and the family of admissible branching mechanisms.
More precisely, denote $\N^{\psi_t}$ the excursion measure induced by $\psi_t$.
Then the process $\{\T_q^t,\; q\in \cT_t\}$ is Markovian  under  $\mN^{\psi_t}$
and  $\N^{\psi_t}(\T_q^t\in d\T)=\N^{\psi_q}(d\T)$. In addition, the special
Markov property holds. Roughly speaking, it gives the conditional distribution
of the tree of individuals with marked ancestors with respect to the tree of
individuals with no marked ancestor.

Due to the consistency property, there is a decreasing tree-valued Markov process
$\{\T_q, q\in\cT\}$ such that $\T_q$ is a $\psi_q$-L\'evy tree.
Let $\bN^\Psi$ be the law of $\{\T_q, q\in\cT\}$. Define the total mass of $\T_q$
and the ascension time  by
$\sigma_q=\bm^{\T_q}(\T_q)$ and $A=\inf\{q\in \cT; \sigma_{q}<+\infty\}.$
The distribution of $\sigma_t$ condition on $\T_q$ and of $A$ under $\bN^\Psi$
are respectively given in Lemmas \ref{lemsigma} and \ref{lemA}.
Conditional on $A$, the distribution of the functionals of the tree at the
ascension time $\T_A$ is considered in Theorem \ref{Theota} and Proposition
\ref{PropAinf}. This generalizes results in \cite{AD12} to admissible
branching mechanisms. An expression for the distribution of the height of the tree
is also given in Proposition \ref{Ah} which is a direct extension of that in \cite{ADH0a}.

Under a regular condition on the admissible family,
 we  prove that the law of the tree-valued process at its ascension time
can be represented in terms of another tree-valued process  obtained by pruning a
critical L\'evy tree conditioned to be infinite;
see Theorem \ref{therep} and Corollary \ref{correp}.

We remark that all the results in this paper are stated using the framework of
real trees but not exploration processes. However  the proof of Theorem \ref{MainA}
relies on explosion process, since Theorem 0.1 and 3.2
of \cite{ADV10} are used explicitly there.

 Let us mention that the study of theory of continuum random trees(CRT) was
 initiated by Aldous \cite{A91,A93}. L\'evy trees, also known as L\'evy CRT,
 were first studied by Le Gall and Le Jan \cite{LL98a,LL98b}, where it was
 shown that L\'evy trees code the genealogy of continuous state branching
 processes (CSBP). Later, in \cite{DL02}, it was shown that Galton-Watson
 trees which code the genealogy of Galton-Watson processes,  suitably rescaled,
 converge to L\'evy trees, as rescaled Galton-Watson processes converge to CSBP.
 Then based on \cite{Li14} and the present work,  one may expect to introduce
 the notation ``admissible family'' to study the Galton-Watson processes and
 Galton-Watson trees. A general pruning procedure on Galton-Watson trees may
 be developed,  which is possibly a combination of Aldous and Pitman's pruning
 procedure in \cite{AP98} and Abraham et al.'s pruning procedure in \cite{ADH12a}.
 This gives the third motivation of the present work.
 We will explore these questions in the future.

  The rest of the paper is organized  as follows. In Section \ref{sec:AFBM},
  we introduce and study the admissible family of branching mechanisms.
  We recall some notations and results on real trees and L\'evy trees
  in Section \ref{sec:levy}. Based on the study of admissible family, in
  Section \ref{sec:prune}, the pruning procedure will be given and the
  marginal distributions of the pruning  process are studied. The evolution
  of the tree-valued process will be explored in Section \ref{sec:tree process}.
  Finally, in the last section,  we construct a tree-valued process by pruning
  a critical L\'evy tree conditioned to be infinite and  get the representation
  of the tree at the  ascension time.

\section{ Admissible family of branching mechanisms} \label{sec:AFBM}
Throughout the paper, for $-\infty\leq a\leq b\leq +\infty$, we make the convention
$\int_a^b=\int_{(a,b)}.$

The admissible family of branching mechanisms was first introduced by Li \cite{Li14}.
Suppose that $\mathfrak{T}\subset{\mbb R}$ is an interval and
$\Psi=\{\psi_q, q\in \cT\}$ is a family of branching mechanisms,
where $\psi_q$ is given by
\[
\psi_q(\lz)=b_q\lz+c\lz^2+\int_0^{\infty}(\expp{-\lz z}-1+\lz z)m_q(dz),\quad\lz\geq0
\]
 with parameters $(b, m)=(b_q, m_q)$ for $q\in \cT$ such that $b_q\in\mR$
 and $\int (z\wedge z^2) m_q(dz)<\infty$.
  \begin{defi}\label{defbranching} [Li(2014)]
  We call $\{\psi_q,q\in \cT\}$ an \textit{admissible family} if for each $\lz>0$ the
  function $q\mapsto\psi_q(\lz)$ is increasing and continuously differentiable with
\beqlb\label{zeta}
\zeta_q(\lz):=\frac{\partial}{\partial q}\psi_q(\lz)
=\beta_q\lz+\int_0^{\infty}(1-\expp{-z\lz})n_q(dz), \quad q\in \cT,\quad \lz>0,
\eeqlb
where $\beta_q\geq0$ and $(1\wedge z)n_q(dz)$ is a finite kernel
from $\cT$ to $(0, \infty)$ satisfying
\beqlb\label{betan1}
\int_t^q \beta_{\theta}d\theta
+\int_t^qd\theta\int_0^{\infty}zn_{\theta}(dz)<\infty,\quad q\geq t\in \cT.
\eeqlb
\end{defi}
\begin{remark}In fact,  it is assumed in \cite{Li14} that $q\mapsto\psi_q(\lz)$
is decreasing and $\zeta_q(\lz)=-\frac{\partial}{\partial q}\psi_q(\lz).$
In that case, we will get an increasing tree-valued process.
\end{remark}
\begin{remark}For the purpose of this work, we also weaken the assumptions
on $\beta_{q}$ and $n_q(dz)$. In \cite{Li14}, it is assumed that
\beqlb\label{betan2}
\sup_{t\leq \theta\leq q}\l(\beta_{\theta}+\int_0^{\infty}zn_{\theta}(dz)\r)
<\infty,\quad q\geq t\in \cT,
\eeqlb
which is essential there. If we assume (\ref{betan2}), some interesting cases
of pruning L\'evy trees may be excluded. See Example \ref{exaAD12} below as an example.
\end{remark}

\begin{remark} It is also possible to assume that
 \[
 \psi_q(\lz)=b_q\lz+c\lz^2+\int_0^{\infty}(\expp{-\lz z}-1+\lz z\ind_{\{z\leq1\}})m_q(dz)
 \]
 with parameters $(b, m)=(b_q, m_q)$ for $q\in \cT$ such that $b_q\in\mR$ and
 $\int (1\wedge z^2) m_q(dz)<\infty$. Then (\ref{betan1}) would be replaced by
 \[
 \int_t^q \beta_{\theta}d\theta+\int_t^qd\theta\int_{(0,1]}zn_{\theta}(dz)<\infty,
 \quad q\geq t\in \cT.
 \]
 We assume further that $\psi_q$ is conservative; i.e.
 $\int_{(0,\ez]}\frac{d\lz}{|\psi_q(\lz)|}=+\infty$ for all $\ez>0$.
 We conjecture that all results in this work can be deduced in this framework.
\end{remark}

\medskip

In the following we  give some examples of the admissible family of branching mechanisms.
\begin{example}\label{exaAD12}
Let $\psi$ be defined in (\ref{eq:psi}). Abraham and Delmas \cite{AD12} considered
$\psi_{q}(\lz)=\psi(q+\lz)-\psi(q), q\in \Theta^{\psi}$,
where $\Theta^{\psi}$ is the set of $\theta\in {\mbb R}$ such that
$\int_1^{\infty}\expp{-\theta z}m(dz)<\infty$. Then $\{\psi_{q}, q\in \Theta^{\psi}\}$
is an admissible family with
\[
b_q=b+2cq+\int_0^{\infty}z(1-\expp{-zq})m(dz),\quad
m_q(dz)=\expp{-zq}m(dz),
\]
and
\[
\beta_q=2c,\quad n_{q}(dz)=z\expp{-zq}m(dz).
\]
Note that $\Theta^{\psi}=[\theta_{\infty}, +\infty)$ or $(\theta_{\infty}, +\infty)$
for some $\theta_{\infty}\in[-\infty,0].$
However, in the case of $\Theta^{\psi}=[\theta_{\infty}, +\infty) $, $n_{\theta_{\infty}}(dz)$
may fail to satisfy (\ref{betan2}). A sufficient condition that (\ref{betan2}) holds is
$\int_1^{\infty}z^2\expp{-\theta_{\infty} z}m(dz)<\infty.$  We remark here that for
the study of the ascension process, we always exclude the case of
$\Theta^{\psi}=[\theta_{\infty}, +\infty)$;
see Remark \ref{tinfty} in Section \ref{sec:tree process} below.
\end{example}

\begin{example}
Let $\psi$ be defined in (\ref{eq:psi}). Let $f\geq0$ be a bounded decreasing function on
 $\mR$ with bounded derivative and $\sup_{x\geq0}|xf'(x)|<+\infty$.
Let $g$ be a  differentiable increasing function on $\mR$. For $q\in\mR$,
let $\psi_q$ be a branching mechanism with parameters $(b_q, m_q)$ defined by
\[
b_q=b+g(q)+\int_0^{\infty}(f(0)-f(zq))zm(dz),\quad
m_q(dz)=f(qz)m(dz).
\]
Then one can check that $\{\psi_q, q\in \mR\}$ is an admissible family of branching mechanisms with
\[
\frac{\partial}{\partial q}\psi_q(\lz)
=g'(q)\lz-\int_0^{\infty}(1-\expp{-z\lz})zf'(qz)m(dz),\quad q\in\mR,\; \lz\geq0,
\]
  and
  \[
  \beta_q=g'(q),\quad n_q(dz)=-zf'(qz)m(dz).
  \]
 In particular, if $m=0$, then $\psi_q(\lz)=(b+g(q))\lz+c\lz^2.$
 If $f\equiv1$, then $\psi_q(\lz)=\psi(\lz)+g(q)\lz.$
\end{example}


\begin{example}\label{exanobody}
Let $\cT_-\subset(-\infty, 0]$ be an interval and  let $\{\psi_q, q\in\cT_-\}$
be an admissible family of branching mechanisms with  parameters $(b_q, m_q)$.
Assume that $0\in\cT_-$ and $\psi_0$ is critical. Note $\eta_q$ the largest root of $\psi_q(s)=0$.
For $q\in -\cT_-:=\{-t, t\in\cT_-\}$, define $\psi_q(\lz)=\psi_{-q}(\lz+\eta_{-q})$.
Then we have $\{\psi_q, q\in \cT_-\cup(-\cT_-)\}$
is an admissible family of branching mechanisms such that for $q\in -\cT_-$,
\[
b_q=b_{-q}+2c\eta_{-q}+\int_0^{\infty}(1-\expp{-z\eta_{-q}})zm_{-q}(dz),
\quad m_q=\expp{-z\eta_{-q}}m_{-q}(dz).
\]
\end{example}

\medskip

 Now we  show how to get pruning parameters for a given admissible family of branching mechanisms.
  Without loss of generality, we always assume that
 $\psi_t\neq \psi_q$ for $t\neq q\in \cT.$
 It follows from  Definition \ref{defbranching} that for $q\geq t\in \cT,$
\beqlb\label{b_q}
b_q=b_t+\int_t^q \beta_{\theta}d\theta+\int_t^qd\theta\int_0^{\infty}zn_{\theta}(dz)
\eeqlb
and
\beqlb\label{m_q}
m_t(dz)=m_q(dz)+\int_{\{t\leq \theta< q\}}n_{\theta}(dz)d\theta.
\eeqlb
\begin{remark} By (\ref{b_q}) one can see $q\mapsto b_q$ is a continuous
increasing function on $\cT$. In particular, $b_t=b_q$ implies
$\psi_t=\psi_q$ and  vice versa.
\end{remark}
For $t\in \cT$, note $\cT_t=\cT\cap [t,+\infty)$.
Using (\ref{m_q}), we get for any $q\in \cT_t$,
$m_q(dz)\ll m_t(dz)\ \text{on}\ (0,\infty).$
Denote by $m_z(t, q)$ the corresponding Radon-Nikodym derivative; i.e.
$m_q(dz)=m_z({t,q})m_t(dz).$
Then we have
\beqlb\label{m_q1}
\int_{t\leq\theta<q}n_{\theta}(dz)d\theta=\left(1-m_z({t,q})\right)m_t(dz),\quad q\in \cT_t,
\eeqlb
which implies $m_{t}(dz)$-a.e.,
\beqlb\label{m_q2}
m_z({t,q})\leq 1 \text{ and } q\mapsto m_z({t,q}) \text{ is decreasing.}
\eeqlb
Furthermore, we yield for $ t\leq\theta\leq q$, $m_{t}(dz)$-a.e.,
\beqlb
\label{m_q3}
m_z(t, q)=m_z(t, \theta)m_z(\theta,q).
\eeqlb
Equation \reff{m_q2},\reff{m_q3} hold $m_{t}(dz)$-a.e.
Since the uncountable union of the null set is not necessarily a null set, then we make the following assumptions:
 \begin{itemize}
 \item[(H1)]  For every $z\in (0,\infty)$ and $t\in \cT$, $$m_z({t,q})\leq 1
 \text{ and } q\mapsto m_z({t,q}) \text{ is decreasing.}$$

 \item[(H2)] For every $z\in (0,\infty)$ and $ t\leq\theta\leq q\in \cT$,
 $$m_z(t, q)=m_z(t, \theta)m_z(\theta,q).$$

 \item[(H3)] For every $q\in \cT$, $$\int^{\infty}\frac{d\lz}{\psi_q(\lz)}<+\infty.$$

 \end{itemize}

\noindent By (H1), we can define a measure $m_z(t, dq)$  on $\cT_t$  by
\beqlb\label{prunmeasure}m_z(t, [t,q])=-\ln (m_z(t, q)),\eeqlb
which induces the pruning measure on branching nodes of infinite degree of a
$\psi_t$-L\'evy tree. By (H2) we have a tree-valued Markov processes on $\cT$.
(H3) is used to ensure that all trees are locally compact.

\medskip

\begin{center}
{\it From now on, we assume that (H1-3) are in force.}
\end{center}

\medskip

\section{Real trees and L\'evy trees}
\label{sec:levy}
In the section, we recall some basic notations and facts on real trees and L\'evy trees.
We mainly follow from Section 2 in \cite{AD13} or \cite{ADH0a}.
\subsection{Notations}

Let $(E,d)$ be a metric  Polish space.  We denote by $M_{f}(E)$ (resp.
$M_f^{\text{loc}}(E)$) the space of all finite (resp.  locally finite)
Borel measures  on $E$.  For $x\in E$,  let $\delta_x$ denote  the Dirac
measure  at point  $x$.  For $\mu\in M_f^{\text{loc}}(E)$  and $f$  a
non-negative measurable function, we set $\langle \mu,f \rangle =\int
f(x) \, \mu(dx)= \mu(f)$.

\subsection{Real trees}

We refer to \cite{E08} or  \cite{Le05} for a general presentation of
random  real trees. A metric
space $(\T,d)$ is a {\it real tree} if the following properties are
satisfied: for every $s,t\in \T$,
\begin{enumerate}
	\item[(i)]  there is a unique isometric map $f_{s,t}:[0,d(s,t)]\rightarrow\T$
	such that $f_{s,t}(0)=s$ and $f_{s,t}(d(s,t))$ $=t$.
	\item[(ii)] if $q$ is a continuous injective map,
$q:[0,1]\rightarrow \T$ such that $q(0)=s$ and $q(1)=t$, then
$q([0,1])=f_{s,t}([0,d(s,t)])$.
\end{enumerate}
If $s,t\in \T$, we will note $\llbracket s,t \rrbracket$
the range of the
isometric map  $f_{s,t}$  and  $\llbracket s,t
\llbracket$  for $\llbracket s,t \rrbracket \backslash \{t\}$.

 $(\T,d,\emptyset)$ is called a rooted real tree with root
$\emptyset$ if $(\T,d)$ is a real tree and $\emptyset\in \T$ is a
distinguished vertex. For every $x\in \T$, $\lb\emptyset ,x\rb$ is
interpreted as the ancestral line of vertex $x$.
 The degree $n(x)$ is  the number of  connected components of $\T\setminus\{x\}$  and the
number of children of $x\neq  \emptyset$ is $\kappa_x=n(x)-1$ and of the
root is  $\kappa_\emptyset=n(\emptyset)$.
 Note $ \Lf (\T)=\{x\in \T,\,  \kappa_x=0\}$, $\Br(\T)=\{x\in \T, \, \kappa_x\geq 2\}$,
$\Br_\infty(\T) =  \{ x\in \T,\, \kappa_x = \infty \}  $ respectively the set of leaves,
branching points and  infinite branching points.
 The skeleton of $\T$ is the  set of points in the
tree that aren't leaves: $\Sk(\T)=\T\backslash \Lf (\T)$.  The trace of the
Borel $\sigma$-field of  $\T$ restricted to $\Sk(\T)$ is  generated by the
sets  $\llbracket  s,s'  \rrbracket$;  $s,s' \in  \Sk(\T)$.   One
defines  uniquely a  $\sigma$-finite  Borel measure  $\ell^{\T}$ on  $\T$,
called the length measure of $\T$, such that
$\ell^{\T}(\Lf (\T)) = 0$ and
$ \ell^{\T}(\llbracket s,s' \rrbracket)=d(s,s').$


\subsection{Measured rooted   real trees}
  We briefly review the
Gromov-Hausdorff-Prohorov  metric on  rooted measured metric space presented in
 \cite{ADH0b}; see also \cite{DW07} and \cite{GPW08} for some related works.

\par Let $(X,d)$ be a Polish metric space.  Recall
$d_\text{H}(A,B)$ the Hausdorff distance between $A$ and $B$ for $A,B\in \cb(X)$
and $d_\text{P}(\mu,\nu)$
the Prohorov distance between $\mu$ and $\nu$ for $\mu,\nu \in M_f(X)$.
A  rooted measured  metric space $\cx = (X, d,  \emptyset,\mu)$ is a
metric space $(X, d)$ with a distinguished element $\emptyset\in X$
and  a locally  finite Borel  measure  $\mu\in M_f^{\text{loc}}(X)$.
Let $\cx$ and $\cx'$ be two compact
 rooted measured metric spaces, and define
\[
 d_{\text{GHP}}^c(\cx,\cx') = \inf_{\Phi,\Phi',Z} \left(d_\text{H}^Z(\Phi(X),\Phi'(X'))
 +d^Z(\Phi(\emptyset),\Phi'(\emptyset')) + d_\text{P}^Z(\Phi_* \mu,\Phi_*'
\mu') \right),
\]
where the infimum is taken over all isometric embedding $\Phi:X\rightarrow
Z$ and $\Phi':X'\rightarrow Z$ into some common Polish metric space
$(Z,d^Z)$ and $\Phi_* \mu$ is the measure $\mu$ transported by $\Phi$.

If $\cx$ is a rooted measured metric space, then for
$r\geq 0$  we will consider  its restriction to  the ball of  radius $r$
centered   at  $\emptyset$,  $\cx^{(r)}=(X^{(r)},   d^{(r)},  \emptyset,
\mu^{(r)})$, where
$X^{(r)}=\{x\in X, d(\emptyset,x)\leq r\}$
with $d^{(r)}$ and $\mu^{(r)}$ defined in an obvious way.

By a measured rooted
 real tree (MRRT) $(\T, d, \emptyset, \bm)$, we mean $(\T,d, \emptyset)$ is a
locally compact  rooted real tree and $  \bm\in\cm_f^\text{loc}(\T)$ is a
locally  finite measure  on $\T$. When there is no confusion, we will simply
write  $\T$  for  $(\T, d,  \emptyset,  \bm)$.
We define   for two MRRTs $\T_1, \T_2$ the Gromov-Hausdorff-Prohorov (GHP)
 metric as follows:
\[
 d_{\text{GHP}}(\T_1,\T_2) = \int_0^\infty \expp{-r} \left(1 \wedge
d^c_{\text{GHP}}\left(\T_1^{(r)},\T_2^{(r)}\right)
\right) \ dr.
\]
 $\T_1$ and
$\T_2$  are  said  GHP-isometric if $ d_{\text{GHP}}(\T_1,\T_2)=0.$
Denote by  $\mT$ the set  of (GHP-isometry classes  of) MRRTs
$(\T, d, \emptyset, \bm)$. According   to  Corollary  2.8   in
\cite{ADH0b}, $(\mT, d_{\text{GHP}})$ is a Polish metric space.

\subsection{Grafting procedure}
\label{sec:graft}
Let $\T$ be  a measured  rooted real  tree and let
$((\T_i,x_i),i\in  I)$ be  a finite  or countable  family of  elements of
$\mT\times \T$.   We define the real  tree obtained by  grafting the trees
$\T_i$  on $\T$  at  point  $x_i$. We  set  $\hat {\T}  =  \T \sqcup  \left(
  \bigsqcup_{i\in  I} \T_i\backslash\{\emptyset^{\T_i}\}  \right)  $ where
the symbol $\sqcup$  means that we choose for  the sets $(\T_i)_{i\in I}$
representatives  of  GHP-isometry classes  in  $\mT$  which are  disjoint
subsets of some common set and that we perform the disjoint union of all
these sets. We set $\emptyset^{\hat{\T}}=\emptyset^{\T}$. The set $\hat{\T}$ is
endowed with the following metric $d^{\hat{\T}}$: if $s,t\in \hat{\T}$,
\begin{equation*}
d^{\hat{\T}} (s,t) =
\begin{cases}
d^{\T}(s,t)\ & \text{if}\ s,t\in \T, \\
d^{\T}(s,x_i)+d^{{\T}_i}(\emptyset^{{\T}_i},t)\ & \text{if}\
s\in {\T},\ t\in {\T}_i\backslash\{\emptyset^{{\T}_i}\}, \\
d^{{\T}_i}(s,t)\ & \text{if}\ s,t\in T_i\backslash\{\emptyset^{{\T}_i}\},\\
d^{\T}(x_i,x_j)+d^{{\T}_j}(\emptyset^{{\T}_j},s)+d^{{\T}_i}
(\emptyset^{{\T}_i},t)\
& \text{if}\ i\neq j \ \text{and}\ s\in {\T}_j\backslash\{\emptyset^{{\T}_j}\},
\ t\in {\T}_i\backslash\{\emptyset^{{\T}_i}\}.
\end{cases}
\end{equation*}
We define the mass measure on $\hat{\T}$ by
\[
\mathbf{m}^{\hat{\T}}=\mathbf{m}^{\T}+\sum_{i\in I}\left(\ind_{
 {\T}_i\backslash\{\emptyset^{{\T}_i}\}} \mathbf{m}^{{\T}_i}+
\mathbf{m}^{{\T}_i}(\{\emptyset^{{\T}_i}\}) \delta_{x_i}\right).
\]
Then $(\hat {{\T}},d^{\hat{\T}},\emptyset^{\hat{\T}})$ is still a complete rooted
real tree (Notice that it is not always true that $\hat {T}$ remains locally
compact  or that $\mathbf{m}^{\hat{\T}}$ is a
locally finite measure on $\hat{\T}$). We use
$\T \otimes_{i\in I}(\T_i,x_i) = (\hat{\T},d^{\hat{\T}},\emptyset^{\hat{\T}},
\mathbf{m}^{\hat{\T}} ) $ for the grafted tree
with the convention that ${\T}\otimes_{i\in I}({\T}_i,x_i)=\T$
for $I=\emptyset$. If $\varphi$ is an isometry from ${\T}$ onto ${\T}'$, then
${\T} \otimes_{i\in I}({\T}_i,x_i) $ and ${\T}' \otimes_{i\in
  I}({\T}_i,\varphi(x_i)) $ are also isometric. Therefore, the grafting
procedure is well defined on $\mT$.

\subsection{Sub-trees above a given level}
 For $\T\in \mT$, define $H_{\text{max}}(\T)=\sup_{x\in \T}
 d^\T(\emptyset^\T,x)$ the height of $\T$ and for
$a\geq 0$,
\[
\T^{(a)}=\{x\in \T, \, d(\emptyset,x)\le a\}
\quad\mbox{and}\quad \T(a)=\{x\in \T, \, d(\emptyset,x)=a\}
\]
the restriction of the tree $\T$ under level $a$ and the set of
vertices of $\T$ at level $a$ respectively.   We denote by  $(\T^{i,\circ},i\in I)$
the connected components of $\T\setminus  \T^{(a)}$.
Let  $\emptyset_i$  be  the most recent common ancestor  of  all  the
vertices    of    $\T^{i,\circ}$.    We    consider    the   real    tree
$\T^i=\T^{i,\circ}\cup\{\emptyset _i\}$ rooted at point $\emptyset_i$ with
mass  measure  $\bm^{\T^i}$ defined  as  the  restriction  of $\bm^\T$  to
$\T^{i,\circ}$ and $\bm^{\T^i}(\emptyset_i)=0$. Notice that
$\T=\T^{(a)} \circledast_{i\in I}(\T^i,\emptyset_i) $.
We will consider the point measure on $\T\times \mT$:
\begin{equation}
   \label{eq:def-cna}
\cn_a^\T=\sum_{i\in I}\delta_{(\emptyset_i,\T^i)}.
\end{equation}

\subsection{Excursion measure of a L\'evy tree}
\label{sec:exc-meas}
Recall (\ref{eq:psi}). We say $\psi$ is subcritical, critical or super-critical
if $b>0$, $b=0$ or $b<0$, respectively.
 In particular, we say $\psi$ is (sub)critical, if $b\geq0$. We assume the
 Grey condition holds:
\begin{equation}\label{Grey}
\int^{+\infty}\frac{d\lambda}{\psi(\lambda)}<+\infty.
\end{equation}

\begin{remark}
The  Grey  condition   is  used to
ensure that the  corresponding L\'{e}vy tree is locally compact and also
implies $c>0$ or $\int_{(0,1)} \ell m(d\ell)=+\infty $
which is equivalent to the fact that the L\'{e}vy process with index
$\psi$ is of infinite variation.
\end{remark}

Let $v^{\psi}$ be
the unique non-negative solution of the equation
\beqlb\label{va}
\int_{v^{\psi}(a)}^{+\infty}\frac{d\lambda}{\psi(\lambda)}=a.
\eeqlb

Results  from  \cite{DL05} in  the  (sub)critical  cases, where height
functions are introduced to  code the compact real trees, can be extended to the
super-critical cases; see \cite{ADH0a}. We  recall the results as follows.
 Note $\N^\psi[d\ct]$
on $\mT$ the excursion  measure of a  L\'{e}vy tree. A $\psi$-L\'evy tree is a
``random'' tree  with law $\mN^{\psi}$  and  the following
properties:
\begin{enumerate}
   \item[(i)] \textbf{Height.} For all $a>0$,
     $\N^\psi[H_{\text{max}}(\ct)>a]=v^{\psi}(a)$.

\item[(ii)] \textbf{Mass measure.} The mass measure $\bm^\ct$ is supported on
$\Lf (\ct)$, $\N^\psi[d\ct]$-a.e.

   \item[(iii)] \textbf{Local time.}
There exists a  process $\{\ell^a, a\geq 0\}$ supported on
$\ct(a)$ which is c\`adl\`ag for the weak topology on the set of finite measures
on $\T$ such that $\N^\psi[d\ct]$-a.e.
$\mathbf{m}^{\ct}(dx) = \int_0^\infty \ell^a(dx) \, da,$
and $\ell^0=0$, $\inf\{a > 0 ; \ell^a = 0\}=
\sup\{a \geq 0 ; \ell^a\neq 0\}=H_{\text{max}}(\ct)$.
\item[(iv)]  \textbf{Branching property.}
Given $\ct^{(a)}$ for any $a>0$, $\cn_a^{\ct}(dx,d\ct')$ under
$\N^\psi[d\ct|H_{\text{max}}(\ct)>a]$ is a Poisson
point measure on $\T(a)\times \mT$ with intensity
$\ell^a(dx)\N^\psi[d\T']$.
\item[(v)] \textbf{Branching points.} $\N^\psi[d\T]$-a.e.,
the branching points of $\ct$ have 2 children or an
  infinite number of children.
\begin{itemize}
\item The set of binary branching points (i.e. with 2 children) is empty
  $\N^\psi$-a.e. if $c=0$ and is a countable dense subset of $\ct$
  if $c>0$;
\item The set $\Br_\infty(\ct)$  of infinite branching points is
  nonempty with $\N^\psi$-positive measure if and only if $m\ne 0$. If
  $\langle  m,1\rangle=+\infty$, the set  $\Br_\infty(\ct)$ is
  $\N^\psi$-a.e. a countable dense subset of $\ct$.
\end{itemize}
\item[(vi)] \textbf{Mass of the nodes.}
The set $\{d(\emptyset,x),\ x\in \Br_\infty(\ct)\}$ coincides
$\N^\psi$-a.e. with the set of discontinuity times of the mapping $a\mapsto
\ell^a$. Moreover, $\N^\psi$-a.e., for every such discontinuity time $b$, there
is a unique $x_b\in \Br_\infty(\ct)$ with $d(\emptyset, x_b)=b$ and
$\Delta_b>0$, such that
$\ell^b = \ell^{b-} + \Delta_b \delta_{x_b}.$
\end{enumerate}

In order to stress the dependence on $\ct$, we may write $\ell^{a, \ct}$
for $\ell^a$.
We set $\sigma^\ct$  or simply $\sigma$ when there  is no confusion, for
the
total mass of the mass measure on $\ct$:
\begin{equation}
\label{eq:s=la}
\sigma=\bm^{\ct}(\ct).
\end{equation}
Notice that  $\bm^\ct(\{x\})=0$ for
any $x\in \T$.

\subsection{Related measures on L\'evy trees}
\label{sec:meas-LT}

We  define a  probability  measure on  $\mT$  as follows.   Let $r>0$  and
$\sum_{k\in\ck}\delta_{\ct^k}$ be  a Poisson random measure  on $\mT$ with
intensity  $r\N^\psi$.   Consider $\emptyset$  as  the trivial  MRRT reduced
to the  root with  null mass  measure. Note
$\T=\emptyset   \circledast_{k\in  \ck}(  \ct^k,   \emptyset)$.   Using
Property (i) as well as \reff{sigma} below, one easily get that $\ct$
is  a  locally  compact MRRT,  and thus  belongs to
$\mT$. We denote by $\mP^\psi_r$ its distribution. The corresponding local
time, mass  measure and  total    mass are  respectively defined  by
$\ell^a=\sum_{k\in \ck}  \ell^{a, \ct^k}$,  $\bm^\ct=\sum_{k\in  \ck}
\bm^{\ct^k}$ and
$\sigma=\sum_{k\in  \ck}  \sigma^{\ct^k}$.   By  construction,  we  have
$\mP^\psi_r(d\ct)$-a.s.          $\emptyset\in          \Br_\infty(\ct)$,
$\Delta_0=r$   and
$\ell^0=r\delta_\emptyset$.   Under $\mP^\psi_r$  (or  $\N^\psi$), we
define the process $\cz=\{\cz_a,a\geq 0\}$ by
\beqlb\label{Za}
\cz_a=\langle\ell^a,1 \rangle .
\eeqlb
Denote by $\eta$ the largest root of $\psi(s)=0$.
Let $\psi^{-1}:[0, +\infty)\mapsto[\eta, +\infty)$ be the inverse function of $\psi$.
 Notice that (under $\mP_r^\psi$ or $\N^{\psi}$):
$\sigma=\int_0^{+\infty } \cz_a \; da =\bm^\ct(\ct).$
In particular, as $\sigma$ is distributed as the total mass of a CSBP
(accumulated mass over all times) under its canonical measure
(intuitively it describes the distribution of CSBP started at an
infinitesimal mass), we have that
 for $\lz\geq 0$,
\begin{equation}
   \label{sigma}
\N^\psi\left[1-\expp{-\lz \sigma} \right]
=\psi^{-1}(\lz),\quad \N^{\psi}[1-\expp{-\lz\cz_a}]=u^{\psi}(a,\lz),
\end{equation}
where $(u^{\psi}(a, \lz), a\geq0, \lz>0)$ is the unique non-negative solution to
\beqlb\label{ua}
\int_{u^{\psi}(a,\lz)}^{\lz}\frac{dr}{\psi(r)}=a;\quad u^{\psi}(0,\lz)=\lz.
\eeqlb
see e.g. $(29)$ and Lemma $2.4$ in \cite{AD12}.
The semigroup property implies
 \beqlb
\label{uva}u^{\psi}(a, u^{\psi}(a',\lz))
=u^{\psi}(a+a',\lz), \quad \lim_{\lz\rar \infty}u^{\psi}(a,\lz)=v^{\psi}(a).
\eeqlb
Finally, we  recall  the Girsanov transformation in  \cite{AD12}.
 Let  $\theta\in \Theta^\psi$ and $a>0$. We set
\[
M_a^{\psi, \theta}=\exp\left\{\theta \cz_0-
  \theta\cz_a-\psi(\theta)\int_0^a\cz_sds\right\}.
\]
Recall that $\cz_0=\la\ell_0,1\ra=0$ under $\N^\psi$.
For any non-negative measurable functional $F$ defined on $\mT$, we have
for $\theta\in \Theta^\psi$ and $a\geq 0$,
 \begin{equation}
 \label{Gir}
 \E_r^{\psi^\theta}[F(\T^{(a)})]
= \E_r^{\psi} \left[F(\T^{(a)})
M_a^{\psi,\theta} \right]
\quad\text{and}\quad
 \N^{\psi^\theta}[F(\T^{(a)})]
= \N^{\psi} \left[F(\T^{(a)})
M_a^{\psi,\theta}  \right].
 \end{equation}
In particular,
\beqlb
\label{Gir3}
 \N^{\psi^\theta}[F(\T)]
&= \N^{\psi} \left[F(\T)\expp{- \psi(\theta) \sigma}
\ind_{\{\sigma<+\infty \}}\right],
\eeqlb
and by (29) in \cite{ADH0a}, we also have for
 $\theta>0$ such that $\psi(\theta)\geq0$,
\begin{equation}
   \label{Gir2}
 \N^{\psi^{\theta}}\left[1-\exp\left\{\theta\Z_{a}+\psi(\theta)
     \int_0^a\Z_sds\right\}\right]=-\theta.
\end{equation}

\section{A general pruning procedure}
\label{sec:prune}
In this section we define a pruning procedure on a L\'evy tree
associated with the admissible family of branching mechanisms.

Recall (\ref{zeta}), (\ref{prunmeasure}) and $\cT_t=\cT\cap [t,\infty)$
for $t\in \cT$. For $\T\in\mT$, we consider under probability measure $\mQ$
two Poisson random measures $M_t^{ske}(d\theta,dy)$ and $M_t^{nod}(d\theta, dy)$
on  the product space $\cT_t\times\T$ with intensity
\[
\beta_{\theta}d\theta\ell^\T(dy)\quad \text{ and }\quad\sum_{x\in
\text{Br}_\infty(\T)\setminus\{\emptyset\}}m_{\Delta_x} (t, d\theta)\delta_x(dy),
\]
respectively. Then $M_t^{ske}(d\theta,dy)$ characterizes the marks on the skeleton
and $M_t^{nod}(d\theta, dy)$ describes the marks on the nodes of infinite degree.

We define a new Poisson random measure on $\cT_t\times\T$ by
\[
M_t^{\T}(d\theta, dy)=M_t^{ske}(d\theta,dy)+M^{nod}_t(d\theta, dy).
\]

\noindent The pruned tree at time $q$ for $q\in \cT_t$ can thus be defined as
\beqlb\label{T_q}
\T_q^t=\{x\in\ct,\ M_t^\T([t, q]\times \lb\emptyset,x\lb)=0\}
\eeqlb
with  the  induced  metric,   root  $\emptyset$  and  mass  measure restricted to $\T_q^t$.

\begin{remark}
Note that $\mN^{\psi_t}$-a.e. $n(\emptyset)=1$ and $\mP_r^{\psi_t}$-a.s.
$n(\emptyset)=\infty$ with $\Delta_0=r$.
The above definition of $\T_q^t$ indicates that we do not add marks on the root
even though $\emptyset$ is a node of infinite degree with mass $r$.
\end{remark}

For fixed $q\in \cT_t$, $M_t^{\T}([t, q], dy)=M_t^{ske}([t,q],dy)+M_t^{nod}([t, q], dy)$
is also a point measure  on tree $\T$:
\begin{itemize}
\item[(i)] $M_t^{ske}([t,q],dy)$  is a Poisson point measure  with intensity
$\int_t^q\beta_{\theta}d\theta \ell^{\T}(dy)$ on the skeleton of $\T$;

\item[(ii)] The atoms of $M_t^{nod}([t, q], dy)$ give the marked nodes:
each node of infinite degree is marked (or pruned) independently from the
others with probability
\beqlb\label{Pmark}
\mQ\l(M_{t}^{nod}([t, q], \{y\})>0\r)=1-\exp\{-m_{\Delta_y}(t, [t, q])\}=1-m_{\Delta_y}(t,q),
\eeqlb
where $\Delta_y$ is the mass  associated with the node.
\end{itemize}
Thus for fixed $q\in\cT_t$, there exists a measurable functional
$\M_{\alpha_{t,q}, p_{t,q}}$ on $\mT$ such that
\beqlb\label{T_q1}
\T_q^t=\M_{\alpha_{t,q}, p_{t,q}}(\T),
\eeqlb
where $\alpha_{t,q}=\int_t^q\beta_{\theta}d\theta, p_{t,q}=1-m_z(t,q).$
 For $t$ fixed and $q\in\cT_t$, if $\psi_t$ is (sub)critical, then we deduce from Theorems 0.1 and 3.2 in \cite{ADV10} the following result in our setting.
 \begin{lem}\label{lem-adv}
\begin{enumerate}
\item[(i)]
The distribution of $\T_q^t$ under $\mN^{\psi_t}$
is $\mN^{\psi_q}$;

\item[(ii)] Given ${\cal T}\in \mT$, let $\cal M(dx, d{\cal T})=\sum_{i\in I}\dz_{(x_i,  \T_i)}$
be a Poisson random measure on ${\cal T}\times {\mbb T}$ with intensity
$
{\bf m}^{{\cal T}}(dx)\l(\alpha_{t,q}{\mbb N}^{\psi_{t}}[d{\cal T}]
+\int_0^{\infty}p_{t,q}{\mbb P}_z^{\psi_{t}}(d{\cal T})m_t(dz)\r).
$
Then under $\N^{\psi_q}$,
$\T\otimes_{i\in I}(\T_i, x_i)$ has the same distribution as $\T$ under $\N^{\psi_t}$.
\end{enumerate}
\end{lem}
We remark that in the special setting of \cite{AD12} the pruning procedure performs with $\beta_q$
 a positive constant and $m_{\Delta_y}(0,q)=\expp{-q\Delta_y}$.
Our main result in this section is the following theorem which is a generalization of the above
result to the supercritical case.

\begin{theo}\label{MainA}
Assume that $\{\psi_t,t\in \cT\}$ is an admissible family satisfying (H1-3). Then we have
\begin{itemize}
\item[(a)]
The  tree-valued process $\{\T_q^t,\; q\in \cT_t\}$ is a
Markov  process  under  $\mN^{\psi_t}$;

\item[(b)]
For fixed $q\in \cT_t$, the distribution of $\T_q^t$ under $\mN^{\psi_t}$
is $\mN^{\psi_q}$;

\item[(c)] Given ${\cal T}\in \mT$, let $\cal M(dx, d{\cal T})=\sum_{i\in I}\dz_{(x_i,  \T_i)}$
be a Poisson random measure on ${\cal T}\times {\mbb T}$ with intensity
\beqlb\label{Mintensity}
{\bf m}^{{\cal T}}(dx)\l(\int_t^q\beta_{\theta}d\theta{\mbb N}^{\psi_{t}}[d{\cal T}]
+\int_t^qd\theta\int_0^{\infty}n_{\theta}(dz){\mbb P}_z^{\psi_{t}}(d{\cal T})\r).
\eeqlb
Then for $q\in\cT_t$,
$(\T, \M_{\alpha_{t,q}, p_{t,q}}(\T))$ under $\mN^{\psi_t}$ has the same distribution as
$(\tilde{\T}, {\T})$ under $\mN^{\psi_q}$,
where \beqlb\label{Ttilde}\tilde{\T}=\T\otimes_{i\in I}(\T_i, x_i).\eeqlb
\end{itemize}
\end{theo}
\begin{remark}\label{remincreasing}
(c) in Theorem \ref{MainA} is the so-called special Markov property which describes
the two dimensional distribution of the tree-valued process.
One may follow the proof in Appendix A in \cite{H12} to extend (c) to have pruning
times in (sub)critical cases and then follow the arguments in {\it Step 4} below
to extend the result to super-critical cases.
\end{remark}

\proof  The proof will be divided into five steps: \\
 {\it Step 1:} We prove (a). It suffices to study the behavior of $M_t^{\T}(d\theta, dy)$
under $\mQ$. Given a branching node $y\in Br_{\infty}(\T)$, for $t\leq \theta\leq q\in\cT$,
 we have
\beqlb\label{MainA1a}
\ar\ar\mQ\l(M_t^{nod}([\theta, q], \{y\})>0\big{|}M_t^{nod}([t, \theta], \{y\})=0\r)\cr
\ar\ar\quad=\frac{\mQ\l(M_t^{nod}([\theta, q], \{y\})>0,
\, M_t^{nod}([t, \theta], \{y\})=0\r)}{\mQ\l(M_t^{nod}([t, \theta], \{y\})=0\r)}\cr
\ar\ar\quad=\frac{\mQ\l(M_t^{nod}([t, q], \{y\})>0)
-\mQ( M_t^{nod}([t, \theta], \{y\})>0\r)}{\mQ\l(M_t^{nod}([t, \theta], \{y\})=0\r)}\cr
\ar\ar\quad=\frac{m_{\Delta_y}(t, \theta)-m_{\Delta_y}(t,q)}{m_{\Delta_y}(t, \theta)}\cr
\ar\ar\quad=1-m_{\Delta_y}(\theta, q)\cr
\ar\ar\quad=\mQ\l(M_{\theta}^{nod}([\theta, q], \{y\})>0\r),
\eeqlb
where we used (\ref{Pmark}) for the third equality and  assumption (H2) for the fourth.
Similarly, one can prove that for $x\in\T$ and $t\leq \theta\leq q\in\cT$,
\[
\mQ\l(M_t^{ske}([\theta, q], \lb\emptyset, x\lb)>0\big{|}M_t^{ske}([t, \theta], \lb\emptyset, x\lb)=0\r)
=\mQ\l(M_{\theta}^{ske}([\theta, q], \lb\emptyset, x\lb)>0\r).
\]
Then (a) follows readily.

\medskip

 {\it Step 2:}
 It remains to study the super-critical case for the second and the third assertions.
 Without loss of generality, we may assume $t=0$. From now on we shall assume that
 $\psi_0$ is super-critical. For this proof only, we set for
 $x\in \T$, $|x|=d^{\T}(\emptyset, x)$. We also write $\T_q$ for $\T_q^0$.

In this step, however, we consider the desired results for a special subcritical case.
Let $\eta_0>0$ be the maximum root of $\psi_0(s)=0.$ Define
\beqlb\label{psit_0}
\psi_q^{\eta_0}(\lz)=\psi_q(\lz+\eta_0)-\psi_q(\eta_0),\quad \lz\geq0,\quad q\in \cT_0.
\eeqlb
One can check that if $\{\psi_q, q\in\cT_0\}$ is an admissible family satisfying (H1-3),
then $\{\psi_q^{\eta_0}, q\in\cT_0\}$ is also an admissible family with parameters
$((b_q^{\eta_0}, m_q^{\eta_0}), q\in \cT_0)$  satisfying (H1-3) such that
\beqlb\label{bmt_0}
b_q^{\eta_0}=b_q+2c\eta_0+\int_0^{\infty}z(1-\expp{-\eta_0z})m_q(dz),\quad m_q^{\eta_0}(dz)
=\expp{-\eta_0 z}m_q(dz).
\eeqlb
In addition, an application of (\ref{b_q}) and (\ref{m_q}) yields
\beqnn
b_q^{\eta_0}\ar=\ar b_0^{\eta_0}+\int_0^q\beta_{\theta}d\theta
+\int_0^qd{\theta}\int_0^{\infty}zn_{\theta}(dz)
-\int_0^{q}d\theta\int_0^{\infty}z(1-\expp{-\eta_0z})n_{\theta}(dz)\cr
\ar=\ar b_0^{\eta_0}+\int_0^q\beta_{\theta}d\theta
+\int_0^{q}d\theta\int_0^{\infty}z\expp{-\eta_0z}n_{\theta}(dz)\cr
\ar\geq\ar b_0^{\eta_0}.
\eeqnn
Since $b_0^{\eta_0}=\psi'_0(\eta_0)>0$, then $\psi_q^{\eta_0}$ is subcritical
for all $q\in \cT_0$. Moreover, by (\ref{psit_0}) and (\ref{bmt_0}),
\beqlb\label{zeta_0}
\frac{\partial}{\partial q}\psi_q^{\eta_0}(\lz)
\ar=\ar\zeta_q(\lz+\eta_0)-\zeta_q(\eta_0)\cr
\ar=\ar \beta_q\lz+\int_0^{\infty}(1-\expp{-\lz z})\expp{-\eta_0z}n_q(dz),
\eeqlb
and
$$
\frac{m_q^{\eta_0}(dz)}{m_0^{\eta_0}(dz)}=\frac{m_q(dz)}{m_0(dz)}=m_z(0,q),
$$
which implies $\{\psi_q, q\in\cT_0\}$ and $\{\psi_q^{\eta_0}, q\in\cT_0\}$
induce the same pruning parameters $\beta_q$ and $m_z(0,q)$.
Therefore, using Lemma \ref{lem-adv} for (sub)critical case, we get
\begin{itemize}
\item[(b')]$\T_q=\M_{\alpha_{0,q}, p_{0,q}}(\T)$ is a $\psi_q^{\eta_0}$-L\'evy
tree under $\mN^{\psi_0^{\eta_0}}$;

 \item[(c')] Given ${\cal T}\in \mT$, let $\cal M^{\eta_0}(dx,d{\cal T})
 =\sum_{i\in I_{\eta_0}}\dz_{(x_i,  \T_i)}$ be a Poisson point measure on
  ${\cal T}\times {\mbb T}$ with intensity
\beqlb\label{M0intensity}
{\bf m}^{{\cal T}}(dx)\l(\int_0^q\beta_{\theta}d\theta{\mbb N}^{\psi_{0}^{\eta_0}}[d{\cal T}]
+\int_0^qd\theta\int_0^{\infty}\expp{-\eta_0z}n_{\theta}(dz){\mbb P}_z^{\psi_{0}^{\eta_0}}(d{\cal T})\r).
\eeqlb
Then for $q\in\cT_0$,
$(\T, \M_{\alpha_{0,q}, p_{0,q}}(\T))$ under $\mN^{\psi_0^{\eta_0}}$ has the same distribution as
$(\hat{\T}_{\eta_0}, {\T})$ under $\mN^{\psi_q^{\eta_0}}$,
where \beqlb\label{Ttilde0}\hat{\T}_{\eta_0}=\T\otimes_{i\in I_{\eta_0}}(\T_i, x_i).\eeqlb
\end{itemize}

{\it Step 3:} We shall prove (b) when $\psi_0$ is super-critical.
Recall $\T^{(a)}=\{x\in \T; d^{\T}(\emptyset, x)\leq a\}$.
By Girsanov transformation (\ref{Gir}), for any nonnegative function $F$ on $\mT$, we have
\beqlb\label{MainA01}
\mN^{\psi_0}[F(\T_q^{(a)})]&=&\mN^{\psi_0}[F(\M_{\alpha_{0,q}, p_{0,q}}(\T^{(a)}))]
\cr
&=&
\mN^{\psi_0^{\eta_0}}[\expp{\eta_0\Z_a}F(\M_{\alpha_{0,q}, p_{0,q}}(\T^{(a)}))]
\cr
&=&
\mN^{\psi_q^{\eta_0}}[\expp{\eta_0\hat{\Z}_a}F(\T^{(a)})],
\eeqlb
where the last equality  follows from Special Markov property (c') and
\[
\hat{\Z}_a=\la\ell^{a,\, \hat{\T}},1\ra
=\Z_a+\sum_{i\in I_{\eta_0}}\ind_{\{|x_i|\leq a\}}\Z_{a-|x_i|}^{\T_i}
\]
with $\Z_a^{\T_i}=\la \ell^{a, \T_i}, 1\ra.$ Then by the property of Poisson random measure,
\beqlb\label{MainA02}
\mN^{\psi_q^{\eta_0}}[\expp{\eta_0 \hat{\Z}_a}F(\T^{(a)})]
=\mN^{\psi_q^{\eta_0}}[\expp{\eta_0 \Z_a}F(\T^{(a)})H(a,\eta_0)]
\eeqlb
with
\beqnn
H(a, \eta_0)\ar=\ar \mN^{\psi_q^{\eta_0}}\l
[\expp{\eta_0\sum_{i\in I_{\eta_0}}\ind_{\{|x_i|\leq a\}}\Z_{a-|x_i|}^{\T_i}}\bigg{|}\T\r]\\
\ar=\ar\exp\bigg{\{}-\int_{\T^{(a)}}{\bf m}^{\T}(dx)\bigg{(}
\int_0^q\beta_{\theta}d\theta\mN^{\psi_0^{\eta_0}}\l[1-\expp{\eta_0\Z_{a-|x|}}\r]
\cr\ar\ar
\qquad\qquad+\int_0^{q}d{\theta}\int_0^{\infty}\expp{-\eta_0z}n_{\theta}(dz)
\l(1-\expp{-z\mN^{\psi_0^{\eta_0}}[1-\expp{\eta_0\Z_{a-|x|}}]}\r)\bigg{)}\bigg{\}}.
\eeqnn
Thanks to (\ref{Gir2}) and the fact that $\psi_0(\eta_0)=0$, we get
\beqlb\label{Nt0}
\mN^{\psi_0^{\eta_0}}\l[1-\expp{\eta_0\Z_{a-|x|}}\r]=-\eta_0,
\eeqlb
which implies
\beqnn
H(a, \eta_0)=\exp\bigg{\{}-\int_{\T^{(a)}}{\bf m}^{\T}(dx)\bigg{(}
-\int_0^q\beta_{\theta}d\theta \eta_0
+\int_0^{q}d{\theta}\int_0^{\infty}\expp{-\eta_0z}n_{\theta}(dz)(1-\expp{z\eta_0})\bigg{)}\bigg{\}}.
\eeqnn
Since $\bm^{\T}(\T^{(a)})=\int_0^{a}\Z_sds$ and
\beqlb\label{MainApsiq}
\psi_q(\eta_0)\ar=\ar\psi_q(\eta_0)-\psi_0(\eta_0)\cr
\ar=\ar \int_0^q\beta_{\theta}d\theta \eta_0
-\int_0^{q}d{\theta}\int_0^{\infty}\expp{-\eta_0z}n_{\theta}(dz)(1-\expp{z\eta_0}),
\eeqlb
then we have
\[
H(a, \eta_0)=\exp\l\{\psi_q(\eta_0){\bf m}^{\T}(\T^{(a)})\r\}=\exp\l\{\psi_q(\eta_0)\int_0^a\Z_sds\r\}.
\]
By (\ref{MainA01}), (\ref{MainA02}) and Girsanov transformation (\ref{Gir}) again, we obtain
\beqnn \mN^{\psi_0}[F(\T_q^{(a)})]\ar=\ar \mN^{\psi_q^{\eta_0}}
[\expp{\eta_0\Z_a+\psi_q(\eta_0)\int_0^a\Z_sds}F(\T^{(a)})]
\cr\ar=\ar\mN^{\psi_q}[F(\T^{(a)})],
\eeqnn
which implies that under $\mN^{\psi_0}$, $\T_q$ is a $\psi_q$-L\'evy tree. We complete the proof of assertion (b).

\medskip

{\it Step 4:} We shall prove (c) when $\psi_0$ is super-critical. Recall (\ref{Ttilde}) and (\ref{Ttilde0}).
Note that
\[
\tilde{\T}^{(a)}=\T^{(a)}\otimes_{i\in I, |x_i|\leq a}(\T_i^{(a-|x_i|)}, x_i).
\]
It suffices  to show that for all $a\geq0$, $(\T^{(a)}, \M_{\alpha_{0,q}, p_{0,q}}(\T^{(a)}))$ under $\mN^{\psi_0}$
has the same distribution as $(\tilde{\T}^{(a)}, \T^{(a)})$ under $\mN^{\psi_q}.$
By (\ref{Gir}), we have for any nonnegative functional $F$ on $\mT^2$,
\beqlb\label{MainA03}
\mN^{\psi_0}\l[F(\T^{(a)}, \M_{\alpha_{0,q}, p_{0,q}}(\T^{(a)}))\r]=
\mN^{\psi_0^{\eta_0}}\l[\expp{\eta_0\Z_a}F(\T^{(a)}, \M_{\alpha_{0,q}, p_{0,q}}(\T^{(a)}))\r].
\eeqlb
By (c') in {\it Step 2}, we deduce that $(\T^{(a)}, \M_{\alpha_{0,q}, p_{0,q}}(\T^{(a)}))$ under $\mN^{\psi_0^{\eta_0}}$
has the same distribution as $(\hat{\T}^{(a)}, \T^{(a)})$ under $\mN^{\psi_q^{\eta_0}}$, where
\[
\hat{\T}^{(a)}=\T^{(a)}\otimes_{i\in I_{\eta_0}, |x_i|\leq a}(\T_i^{(a-|x_i|)}, x_i).
\]
Thus we yield
\beqlb\label{MainA05}
\mN^{\psi_0}\l[F(\T^{(a)}, \M_{\alpha_{0,q}, p_{0,q}}(\T^{(a)}))\r]
=\mN^{\psi_q^{\eta_0}}\l[\expp{\eta_0\hat{\Z}_a}F(\hat{\T}^{(a)}, \T^{(a)})\r].
\eeqlb
We CLAIM that for all $a>0$ and any nonnegative measurable functional $\Phi$ on $\T^{(a)}\times\mT$,
\beqlb\label{MainA04}
\mN^{\psi_q^{\eta_0}}\l[\expp{\eta_0\hat{\Z}_a}F(\T^{(a)})\exp\{-\la\M_a^{\eta_0},\Phi\ra\}\r]
=\mN^{\psi_q}\l[F(\T^{(a)})\exp\{-\la\M_a,\Phi\ra\}\r],\eeqlb
where
\[
\M_a^{\eta_0}(dx, d\T)=\sum_{i\in I_{\eta_0}} \ind_{|x_i|\leq a}\dz_{(x_i,  \T_i^{(a-|x_i|)})}(dx,  d\T),
\]
and
\[
 \M_a(dx,  d\T)=\sum_{i\in I}\ind_{|x_i|\leq a}\dz_{(x_i,  \T_i^{(a-|x_i|)})}(dx,  d\T).
\]
Then we deduce from (\ref{MainA04}) that
\beqnn
\mN^{\psi_q^{\eta_0}}\l[\expp{\eta_0\hat{\Z}_a}F(\hat{\T}^{(a)}, \T^{(a)})\r]\ar
=\ar\mN^{\psi_q^{\eta_0}}\l[\expp{\eta_0\hat{\Z}_a}F(\T^{(a)}
\otimes_{i\in I_{\eta_0}, |x_i|\leq a}(\T_i^{(a-|x_i|)}, x_i), \T^{(a)})\r]\cr
\ar=\ar\mN^{\psi_q}\l[F(\T^{(a)}\otimes_{i\in I, |x_i|\leq a}(\T_i^{(a-|x_i|)}, x_i), \T^{(a)})\r]\cr
\ar=\ar\mN^{\psi_q}\l[F(\tilde{\T}^{(a)}, \T^{(a)})\r],
\eeqnn
which, together with (\ref{MainA05}), gives
\[
\mN^{\psi_0}\l[F(\T^{(a)}, \M_{\alpha_{0,q}, p_{0,q}}(\T^{(a)}))\r]
=\mN^{\psi_q}\l[F(\tilde{\T}^{(a)}, \T^{(a)})\r].
\]
Since $a$ is arbitrary, assertion (c) follows readily.

\bigskip

{\it Step 5:} The remainder of this proof is devoted to (\ref{MainA04}).
Define
\[
g(a, x)=\mN^{\psi_0}\l[1-\expp{-\Phi(x,  \T^{(a-|x|)})}\r].
\]
Then we have
\[
\mP_z^{\psi_0}\l(1-\expp{-\Phi(x, \T^{(a-|x|)})}\r)=1-\expp{-zg(a, x)}.
\]
First, by the property of Poisson random measure, we get
\beqlb\label{MainA04a}
\ar\ar\mN^{\psi_q}\l[F(\T^{(a)})\exp\{-\la\M_a, \Phi\ra\}\r]\cr
\ar\ar\qquad= \mN^{\psi_q}\l[F(\T^{(a)})\exp\l\{-\int_{\T^{(a)}}\bm^{\T^{(a)}}(dx)G(a,x)\r\}\r],
\eeqlb
where
\beqlb\label{MainA04c}
G(a,x)=\l[\int_0^q\beta_{\theta}d\theta g(a, x)
+\int_0^qd\theta\int_0^{\infty}n_{\theta}(dz)\l(1-\expp{-zg(a, x)}\r)\r].
\eeqlb
Then thanks to (\ref{Gir}), we  obtain
\beqlb\label{MainA04d}
g(a, x)
\ar=\ar\mN^{\psi_0}\l[1-\expp{-\Phi(x, \T^{(a-|x|)})}\r]\cr
\ar=\ar\mN^{\psi_0^{\eta_0}}\l[\expp{\eta_0\Z_{a-|x|}}\l(1-\expp{-\Phi(x, \T^{(a-|x|)})}\r)\r]\cr
\ar=\ar\mN^{\psi_0^{\eta_0}}\l[\expp{\eta_0\Z_{a-|x|}}-1+1-\expp{-\Phi(x,  \T^{(a-|x|)})+\eta_0\Z_{a-|x|}}\r]\cr
\ar=\ar \eta_0+\mN^{\psi_0^{\eta_0}}\l[1-\expp{-\Phi(x, \T^{(a-|x|)})+\eta_0\Z_{a-|x|}}\r]=:\eta_0+g_{\eta_0}(a,x),
\eeqlb
where the last equality follows from (\ref{Nt0}). Using (\ref{MainApsiq}) and (\ref{MainA04d}), we have
\beqlb\label{MainA04b}
G(a,x)\ar=\ar\int_0^q\beta_{\theta}d\theta \eta_0
-\int_0^qd\theta\int_0^{\infty}\expp{-z\eta_0}n_{\theta}(dz)\l(1-\expp{z\eta_0}\r)\cr
\ar\ar\quad+\int_0^q\beta_{\theta}d\theta g_{\eta_0}(a,x)
+\int_0^qd\theta\int_0^{\infty}\expp{-z\eta_0}n_{\theta}(dz)\l(1-\expp{-zg_{\eta_0}(a,x)}\r)\cr
\ar=\ar \psi_q(\eta_0)+\int_0^q\beta_{\theta}d\theta g_{\eta_0}(a,x)
+\int_0^qd\theta\int_0^{\infty}\expp{-z\eta_0}n_{\theta}(dz)\l(1-\expp{-zg_{\eta_0}(a,x)}\r)\cr
\ar=:\ar\psi_q(\eta_0)+G_{\eta_0}(a,x).
\eeqlb
Applying (\ref{Gir}) to (\ref{MainA04a}) gives
\beqlb\label{MainA04e}
\ar\ar\mN^{\psi_q}\l[F(\T^{(a)})\exp\{-\la\M_a, \Phi\ra\}\r]\cr
\ar\ar\qquad= \mN^{\psi_q^{\eta_0}}\l[\exp\l\{{\eta_0\Z_{a}}+\psi_q(\eta_0)\int_0^{a}\Z_sds\r\}
F(\T^{(a)})\exp\l\{-\int_{\T^{(a)}}\bm^{\T^{(a)}}(dx)G(a,x)\r\}\r]\cr
\ar\ar\qquad
=\mN^{\psi_q^{\eta_0}}\l[\expp{\eta_0\Z_{a}}
F(\T^{(a)})\exp\l\{-\int_{\T^{(a)}}\bm^{\T^{(a)}}(dx)G_{\eta_0}(a,x)\r\}\r],
\eeqlb
where the last equality follows from (\ref{MainA04b}) and the fact that $\bm^{\T^{(a)}}(\T^{(a)})=\int_0^{a}\Z_sds.$

On the other hand, using the property of Poisson point measure again, we obtain
\beqnn
\ar\ar \mN^{\psi_q^{\eta_0}}\l[\expp{\eta_0\hat{\Z}_a}F(\T^{(a)})\exp\{-\la\Phi, \M_a^{\eta_0}\ra\}\r]\cr
\ar\ar\quad=\mN^{\psi_q^{\eta_0}}\l[\expp{\eta_0{\Z}_a}F(\T^{(a)})
\exp\l\{-\la\Phi, \M_a^{\eta_0}\ra+\eta_0\sum_{i\in I_{\eta_0}}\ind_{|x_i|\leq a}\Z_{a-|x_i|}^{\T_i}\r\}\r]
\cr
\ar\ar\quad=\mN^{\psi_q^{\eta_0}}\l[\expp{\eta_0\Z_{a}}
F(\T^{(a)})\exp\l\{-\int_{\T^{(a)}}\bm^{\T^{(a)}}(dx)G_{\eta_0}(a,x)\r\}\r],
\eeqnn
which, together with (\ref{MainA04e}), implies (\ref{MainA04}).\qed

A direct consequence of Theorem \ref{MainA} is as follows.

\begin{corollary}\label{coromain1}
Assume that $\{\psi_t,t\in \cT\}$ is an admissible family satisfying (H1-3). Then for $r>0$ we have
the  tree-valued process $\{\T_q^t,\; q\in \cT_t\}$ is a
Markov  process  under  $\mP^{\psi_t}_r$ and for fixed $q\in \cT_t$, the distribution of $\T_q^t$ under $\mP_r^{\psi_t}$
is $\mP_r^{\psi_q}$.
\end{corollary}

\section{a tree-valued process}\label{sec:tree process}
From the construction of $\{\T_q^t,q\in\cT_t\}$,  we have $\T^t_{q}\subset\T^t_{p}$
for $ p\leq q\in\cT_t$.   The  process
$\{\T_q^t, q\in\cT_t\}$ is  a  non-increasing process  (for  the inclusion  of
trees). Corollary \ref{coromain1} (rep. Theorem \ref{MainA}) implies that for $t_1\leq t_2\in \cT$,
$\{\T_q^{t_2}, q\in \cT_{t_2}\}$ under $\mP_r^{\psi_{t_2}}$(resp. $\mN^{\psi_{t_2}}$)
has the same distribution as $\{\T_q^{t_1}, q\in \cT_{t_2}\}$
under $\mP_r^{\psi_{t_1}}$(resp. $\mN^{\psi_{t_1}}$). Then there exists a projective limit
$\{\T_t, t\in\cT\}$ which is a tree-valued process such that $\{\T_q, q\in\cT_t\}$
has the same finite dimensional distribution as $\{\T_q^t, q\in \cT_t\}$ under $\mN^{\psi_t}$.
Denote by $\bP^{\Psi}$ and $\bN^{\Psi}$ the distribution and excursion law of $\{\T_t, t\in\cT\}$.
We have for any nonnegative measurable functional $F$,
\[
\bN^{\Psi}[F(\T_q)]=\mN^{\psi_q}[F(\T)].
\]
Set
\[
\sigma_{t}=\bm^{\ct_t}(\ct_t),\quad t\in \cT.
\]
Then one can check that $\{\sigma_t, t\in\cT\}$ is a non-increasing $[0,\infty]$-Markov process.
For $t\in\cT$,
set $\Psi_t=\{\psi_q, q\in\cT_t\}$ and $\Psi_t^{\eta_t}=\{\psi_q^{\eta_t}, q\in\cT_t\}$. We study
in this section the property of the tree-valued process which is direction generalization of Section
6 in \cite{AD12}.
\begin{proposition}\label{coromain2}
 For $t\in \cT$ and any non-negative measurable functional $F$,
\[
\bN^{\Psi_t}[F(\T_{q}, q\in\cT_t)\ind_{\{\sigma_t<\infty\}}]=\bN^{\Psi_t^{\eta_t}}[F(\T_q, q\in\cT_t)].
\]
\end{proposition}
\proof
Recall that  $\{\psi_q, q\in\cT_t\}$ and $\{\psi_q^{\eta_t}, q\in\cT_t\}$ induce the same pruning parameters.
Then the desired result is a direct consequence of the fact
$\mN^{\psi_t}[F(\T)\ind_{\{\sigma<\infty\}}]=\mN^{\psi_t^{\eta_t}}[F(\T)]$; see (\ref{Gir3}).   \qed

 We then study the behavior of $\{\sigma_t, t\in\cT\}.$

\begin{lemma}\label{lemsigma} For $t\leq q\in \cT$ and $\lz\geq0$, we have
\[
\bN^{\Psi}[\expp{-\lz \sigma_t}|\T_{q}]=\exp\{{-\psi_{q}(\psi_t^{-1}(\lz))\sigma_{q}}\}
\]
and $\bN^{\Psi}[\sigma_t<+\infty|\T_{q}]=\exp\{{-\psi_{q}(\psi_t^{-1}(0))\sigma_{q}}\}$.
Moreover, if $\psi_t$ is subcritical, then
\beqlb\label{lemsigma01}
\bN^{\Psi}[\sigma_t|\T_q]=\psi'_q(0)\sigma_q/\psi'_t(0).
\eeqlb
\end{lemma}
\proof Recall (\ref{zeta}) and (\ref{sigma}). Using (c) in Theorem \ref{MainA}, we obtain
\beqnn
\bN^{\Psi}\l[\expp{-\lz \sigma_t}|\T_{q}\r]
\ar=\ar\bN^{\Psi}\l[\expp{-\lz \sigma_{q}-\lz\sum_{i\in I}\sigma_i}|\T_{q}\r]\cr
\ar=\ar \expp{-\lz\sigma_{q}}\expp{- \int_{\T_{q}}\bm^{\T_{q}}(dx)G(\lz)},
\eeqnn
where $\sigma_i=\bm^{\T_i}(\T_i)$ and
\beqnn
G(\lz)\ar=\ar \int_t^{q}\beta_{\theta}d{\theta}\mN^{\psi_t}\l[1-\expp{-\lz\sigma}\r]
+\int_t^{q}d\theta\int_0^\infty n_{\theta}(dz)\mP^{\psi_t}_z(1-\expp{-\lz\sigma})\cr
\ar=\ar \psi_t^{-1}(\lz)\int_t^{q}\beta_{\theta}d{\theta}
+\int_t^{q}d\theta\int_0^\infty n_{\theta}(dz)\l(1-\expp{-z\psi_t^{-1}(\lz)}\r)\cr
\ar=\ar \int_t^q\zeta_{\theta}(\psi_t^{-1}(\lz))d\theta\cr
\ar=\ar \psi_{q}(\psi_t^{-1}(\lz))-\psi_{t}(\psi_t^{-1}(\lz)).
\eeqnn
Hence,
\beqlb\label{lemsigma02}
\bN^{\Psi}[\expp{-\lz \sigma_t}|\T_{q}]
=\exp\{{-\psi_{q}(\psi_t^{-1}(\lz))\sigma_{q}}\}.
\eeqlb
Consequently,
\[
\bN^{\Psi}[\sigma_t<+\infty|\T_{q}]=\lim_{\lz\rar0}\bN^{\Psi}[\expp{-\lz \sigma_t}|\T_{q}]
=\exp\{{-\psi_{q}(\psi_t^{-1}(0))\sigma_{q}}\}.
\]
If $\psi_t$ is subcritical,  then $\bN^{\Psi}$-a.e. $\sigma_t<\infty$, so we conclude that
\[
 \bN^{\Psi}[\sigma_t|\T_q]=\frac{d}{d\lz}\bN^{\Psi}[\expp{-\lz \sigma_t}|\T_{q}]\big{|}_{\lz=0}
 =\psi'_q(0)\sigma_q/\psi_t'(0).
\]
\qed



Recall that $\eta_q$ is the largest root of $\psi_{q}(\lz)=0$. Thus
\[
\eta_q=\lim_{\lz\rar 0+}\psi_q^{-1}(\lz)=\psi_q^{-1}(0).\]
 Define the {\it ascension time}
\[
A=\inf\{t\in \cT; \sigma_{t}<+\infty\}
\]
with the convention that $\inf\{\emptyset\}=\inf\cT=:t_{\infty}.$
Since $\cT$ is an interval, we always assume that
\[
0\in\cT,\quad t_{\infty}<0.
\]
 Recall (\ref{zeta}). Let us consider the following condition:

 \beqlb\label{bmcon1}
 \lim_{t\rar t_{\infty}+}\int^0_{t}\zeta_{\theta}(\lz)d\theta
 =\psi_0(\lz)-\lim_{t\rar t_{\infty}+}\psi_t(\lz)<+\infty,
 \ \text{ for some }\lz>0.
 \eeqlb
 We have
 \begin{proposition} \label{lemtinf}
 $\lim_{q\rar t_{\infty}+}\psi_q^{-1}(0)<\infty$ if and only if (\ref{bmcon1}) holds.
 \end{proposition}
\proof 
{\it ``if'' part:}  Condition (\ref{bmcon1}) implies
\beqlb\label{bmcon2}
\lim_{t\rar t_{\infty}+}\int^0_{t}\beta_{\theta}d\theta
+\lim_{t\rar t_{\infty}+}\int^0_{t}\int_0^{\infty}(1\wedge z)n_{\theta}(dz)d\theta<+\infty.
\eeqlb
By  (\ref{m_q}) , (\ref{bmcon2}) and the fact that $1\wedge z^2\leq 1\wedge z$,
 we have $\sup_{q\in \cT}\int_0^{\infty}(1\wedge z^2)m_q(dz)<\infty$.
Furthermore, thanks to monotonicity of $q\mapsto (1\wedge z^2)m_q(dz)$, there
exists a $\sigma$-finite measure $m_{t_{\infty}}(dz)$ on $(0, \infty)$ such that
\[
\int(1\wedge z^2)m_{t_{\infty}}(dz)<+\infty,
\]
 and  as ${q\rar t_{\infty}}$,
$(1\wedge z^2)m_q(dz){\rar} (1\wedge z^2)m_{t_{\infty}}(dz)$ in $M_f((0,\infty)).$
Therefore, we may define for $\lz\geq 0$,
\[
\psi_{t_{\infty}}(\lz)=b_{t_{\infty}}\lz+c\lz^2+\int_0^{\infty}(\expp{-\lz z}-1
+\lz z\ind_{\{z\leq 1\}})m_{t_{\infty}}(dz)
\]
with
\[
b_{t_{\infty}}
=b_0+\int_1^{\infty}zm_0(dz)-\lim_{q\rar t_{\infty}+}\l(\int^0_{q}\beta_{\theta}d\theta
+\int^0_q\int_0^1z n_{\theta}(dz)d\theta\r).
\]
We deduce that $\psi_{t_{\infty}}(\lz)$ is a convex function. Since $\psi_q$ satisfies (H3)
which implies $c>0$ or $\int_0^1zm_q(dz)=\infty,$ then we have
$\lim_{\lz\rar\infty}\psi_{t_{\infty}}(\lz)=\infty$ and $\psi^{-1}_{t_{\infty}}(0)<\infty$
 is the largest root of $\psi_{t_{\infty}}(\lz)=0$.
Notice that $ \expp{-\lz z}-1+\lz z\ind_{\{z\leq 1\}}\leq 1\wedge z^2$.
 By (\ref{b_q}) and (\ref{m_q}),
\beqnn
\psi_q(\lz)
\ar=\ar
\l(b_0+\int_1^{\infty}zm_0(dz)-\int^0_q\beta_{\theta}d\theta
-\int^0_q\int_0^{1}zn_{\theta}(dz)d\theta\r)\lz\cr
\ar\ar\quad +c\lz^2+\int_0^{\infty}(\expp{-\lz z}-1+\lz z\ind_{\{z\leq 1\}})m_q(dz)\cr
\ar\rar\ar \psi_{t_{\infty}}(\lz),\quad \text{as }q\rar t_{\infty}+.
\eeqnn

\noindent Then we conclude
$\psi_{t_{\infty}}^{-1}(0)=\lim_{q\rar t_{\infty}+}\psi_q^{-1}(0)<\infty.$

\smallskip

{\it ``only if'' part:}  If $\int_{t_{\infty}+}^0\zeta_{\theta}(\lz)d\theta=+\infty$
for some $\lz>0$ (hence for all $\lz>0$), by (\ref{b_q}) and (\ref{m_q}),
  \beqnn
\psi_q(\lz)=\psi_0(\lz)-\int^0_q\zeta_{\theta}(\lz)d\theta
\rar -\infty,\quad \text{as }q\rar t_{\infty}+.
\eeqnn
Then we have
$\lim_{q\rar t_{\infty}+}\psi_q^{-1}(0)=+\infty.$ \qed

Define $\psi_{t_{\infty}}^{-1}(0)=\lim_{q\rar t_{\infty}}\psi_q^{-1}(0)$ and
\beqnn
 \cT_{\infty}=\begin{cases}\cT\cup\{t_{\infty}\},&\psi_{t_{\infty}}^{-1}(0)<+\infty\\
\cT,& \psi_{t_{\infty}}^{-1}(0)=+\infty.
\end{cases}
\eeqnn
\begin{remark}\label{tinfty}From the proof of Proposition \ref{lemtinf} we see that
it is possible to extend the definition of a given admissible family to $\cT_{\infty}$.
For some results in the sequel of this paper, we need to avoid this case by assuming
$t_{\infty}\notin \cT_{\infty}$ (hence $t_{\infty}\notin \cT$).
\end{remark}
Next, we study the distribution of $A$ and $\T_A$.
\begin{lemma}\label{lemA} For $q\in \cT\cup\{t_{\infty}\}$,
\[
\bN^{\Psi}[A>q]=\psi_{q}^{-1}(0),
\]
and
\beqlb\label{lemA1}
\bN^{\Psi}[A=t_{\infty}]=\begin{cases} 0,& t_{\infty}\notin \cT_{\infty}\\
\infty,& t_{\infty}\in\cT_{\infty}.
\end{cases}
\eeqlb
\end{lemma}
\proof Recall (\ref{sigma}). By Lemma \ref{lemsigma}, for $q>t_{\infty}$,
 \beqnn
\bN^{\Psi}[A>q]\ar=\ar \bN^{\Psi}\l[\sigma_{q}=+\infty\r]
\cr\ar=\ar\mN^{\psi_{q}}\l[\sigma=+\infty\r]\cr
\ar=\ar\lim_{\lz\rar0}\mN^{\psi_{q}}\l[1-\expp{-\lz \sigma}\r]\cr
\ar=\ar\lim_{\lz\rar0}\psi_{q}^{-1}(\lz)\cr
\ar=\ar \psi_{q}^{-1}(0).
\eeqnn
Letting $q\rar t_{\infty}$ gives the case of $q=t_{\infty}.$
Using Lemma \ref{lemsigma} again, we obtain
\beqnn
\bN^{\Psi}\l[A=t_{\infty}|\T_0\r]\ar=\ar
\bN^{\Psi}\l[\forall q>t_{\infty}, \sigma_q<+\infty|\T_0\r]\cr
\ar=\ar \lim_{q\rar t_{\infty}}\bN^{\Psi}\l[ \sigma_q<+\infty|\T_0\r]\cr
\ar=\ar \lim_{q\rar t_{\infty}}\expp{-\psi_0(\psi_q^{-1}(0))\sigma_0}\cr
\ar=\ar
\begin{cases}
0,& \text{ if }t_{\infty}\notin \cT_{\infty}\\
\expp{-\sigma_0\psi_0(\psi_{t_{\infty}}^{-1}(0))},& \text{ if }t_{\infty}\in \cT_{\infty}.
\end{cases}
\eeqnn
Notice that $\forall \lz>0, \mN^{\psi_0}[\expp{-\lz \sigma}]=+\infty$,  the desired follows. \qed

\begin{remark}
(\ref{lemA1}) implies that if $t_{\infty}\in \cT_{\infty}$, then
$\bN^{\Psi}[\T_t \text{ is compact for all }t>t_{\infty}]=+\infty.$
If $t_{\infty}\notin\cT_{\infty}$, then $\bN^{\Psi}$-a.e.
there exists $t\in \cT$ such that $\T_q$ is not compact ($\sigma_q=\infty$) for $t>q\in \cT$.
\end{remark}

For the rest of the paper,  we focus on the ascension time $A$ and tree at the ascension time $\T_A$.
Then it is necessary that there exists some point $c\in\cT$, such that for $q<c$, $\psi_q$ is a
supercritical branching mechanism. For this purpose, {\bf from now on, we  always  assume that $\psi_0$ is critical and $t_{\infty}<0$.}
\begin{theo}\label{Theota}
For $q\in (t_{\infty}, 0)$ and any nonnegative measurable functional $F$ on $\mT$,
\beqlb\label{TA}
\bN^{\Psi}[F(\T_A)|A=q]=\psi'_{q}(\eta_{q})
\mN^{\psi_{q}}[F(\T)\sigma\ind_{\{\sigma<\infty\}}]
\eeqlb
and for $\lz\geq0$,
\beqlb\label{sigmaA}
\bN^{\Psi}[\expp{-\lz \sigma_A}|A=q]=\frac{\psi'_{q}(\eta_{q})}{\psi'_{q}(\psi_{q}^{-1}(\lz))}.
\eeqlb
In particular, we have
\beqlb\label{sigmaA1}
\bN^{\Psi}[\sigma_A<\infty|A=q]=1.
\eeqlb
\end{theo}
\proof By Lemma \ref{lemsigma}, we have for every $t_{\infty}<t\leq q<0$,
\beqnn
\bN^{\Psi}\l[F(\T_{q})\ind_{\{A>t\}}\r]
\ar=\ar
\bN^{\Psi}\l[F(\T_{q})\ind_{\{\sigma_t=+\infty\}}\r]
\cr
\ar=\ar \bN^{\Psi}\l[F(\T_{q})\bN^{\Psi}[\sigma_t=+\infty|\T_{q}]\r]\cr
\ar=\ar \bN^{\Psi}\l[F(\T_{q})\l(1-\expp{-\sigma_{q}\psi_{q}(\psi_t^{-1}(0))}\r)\r]\cr
\ar=\ar \bN^{\Psi}\l[F(\T_{q})\l(1-\expp{-\sigma_{q}\psi_{q}(\eta_{t})}\r)\r].
\eeqnn
Since $\eta_t$ is the largest root of $\psi_t(s)=0$, we have the mapping $t\mapsto\eta_{t}$ is differentiable with
\beqlb\label{detat}
\frac{d\eta_{t}}{dt}=-\frac{\zeta_{t}(\eta_{t})}{\psi'_{t}(\eta_{t})}.
\eeqlb
Then we get
\beqnn
\ar\ar\frac{d}{dt}\bN^{\Psi}\l[F(\T_{q})\ind_{\{A>t\}}\r]\cr
\ar\ar\qquad=\bN^{\Psi}
  \l[F(\T_{q})\sigma_{q}\expp{-\sigma_{q}\psi_{q}(\eta_{t})} \r]
  \frac{d\psi_{q}(\eta_{t})}{dt}\cr
\ar\ar\qquad= -\bN^{\Psi}
  \l[F(\T_{q})\sigma_{q}\expp{-\sigma_{q}\psi_{q}(\eta_{t})} \r]
  \frac{\psi'_{q}(\eta_{t})\zeta_{t}(\eta_{t})}{\psi'_{t}(\eta_{t})}.
  \eeqnn
So we have
 \beqnn
 \frac{\bN^{\Psi}\l[F(\T_A), A\in dq\r]}{dq}\ar
 =\ar-\frac{d}{dt}\l(\bN^{\Psi}\l[F(\T_{q})\ind_{\{A>t\}}\r]\r)\bigg{|}_{t=q}\cr
 \ar=\ar\zeta_{t}(\eta_{t})\bN^{\Psi}
  \l[F(\T_{q})\sigma_{q}\ind_{\{\sigma_{q}<+\infty\}}\r].
   \eeqnn
Thus
\beqlb\label{TA1}
\bN^{\Psi}\l[F(\T_A)| A=q\r]=
\frac{ \bN^{\Psi}\l[F(\T_{q})\sigma_{q}\ind_{\{\sigma_{q}<+\infty\}}\r]}
{\bN^{\Psi}\l[\sigma_{q}\ind_{\{\sigma_{q}<+\infty\}}\r]}=
\frac{ \mN^{\psi_{q}}\l[F(\T)\sigma\ind_{\{\sigma<+\infty\}}\r]}
{\mN^{\psi_{q}}\l[\sigma\ind_{\{\sigma<+\infty\}}\r]}.
\eeqlb
Notice that
\[
\mN^{\psi_{q}}\l[\sigma \expp{-r\sigma}\r]=\frac{d}{dr}\mN^{\psi_{q}}\l[1-\expp{-r\sigma}\r]
=\frac{d}{dr}\psi_{q}^{-1}(r)=\frac{1}{\psi'_{q}(\psi_{q}^{-1}(r))}.
\]
Then we have
\[
\mN^{\psi_{q}}\l[\sigma\ind_{\{\sigma<+\infty\}}\r]
=\lim_{r\rar0}\mN^{\psi_{q}}\l[\sigma \expp{-r\sigma}\r]
=\frac{1}{\psi'_{q}(\eta_{q})},
\]
which, together with (\ref{TA1}), implies  (\ref{TA}). Using (\ref{TA1}) again, we deduce
\beqlb\label{sigmaA2}
\bN^{\Psi}\l[\expp{-\lz \sigma_A}|A=q\r]=
\frac{ \mN^{\psi_{q}}\l[\expp{-\lz \sigma}\sigma\r]}
{\mN^{\psi_{q}}\l[\sigma\ind_{\{\sigma<+\infty\}}\r]}
=\frac{\psi'_{q}(\eta_{q})}{\psi'_{q}(\psi_{q}^{-1}(\lz))}.
\eeqlb
 Then (\ref{sigmaA1}) is a direct consequence of  (\ref{sigmaA2}) by letting $\lz\to 0$. \qed

\medskip
\begin{proposition}\label{PropAinf}
Assume $t_{\infty}\in \cT$. Then for any nonnegative measurable functional $F$ on $\mT$,
\[
\bN^{\Psi}\l[F(\T_A)\ind_{\{A=t_{\infty}\}}\r]=\mN^{\psi_{t_{\infty}}^{\eta_{t_{\infty}}}}\l[F(\T)\r],
\]
where $\psi_{t_{\infty}}^{\eta_{t_{\infty}}}(\lz)= \psi_{t_{\infty}}(\eta_{t_{\infty}}+\lz).$
In particular, for $\lz\geq0$,
\[
\bN^{\Psi}\l[(1-\expp{-\lz\sigma_A})\ind_{\{A=t_{\infty}\}}\r] =\psi_{t_{\infty}}^{-1}(\lz)-\eta_{t_{\infty}}.
\]
\end{proposition}
\proof
First we have
\beqlb\label{TA2}
\bN^{\Psi}\l[F(\T_A)\ind_{\{A=t_{\infty}\}}\r]
\ar=\ar \bN^{\Psi}\l[F(\T_{t_{\infty}})\ind_{\{\sigma_{t_{\infty}}<+\infty\}}\r]\cr
\ar=\ar \mN^{\psi_{t_{\infty}}}\l[F(\T)\ind_{\{\sigma<+\infty\}}\r]\cr
\ar=\ar \mN^{\psi_{t_{\infty}}^{\eta_{t_{\infty}}}}\l[F(\T)\r],
\eeqlb where the last equality follows from (\ref{Gir3}).
 Then we deduce from (\ref{TA2}) that
\beqlb\label{TA3}
\bN^{\Psi}\l[(1-\expp{-\lz\sigma_A})\ind_{\{A=t_{\infty}\}}\r]
\ar=\ar \mN^{\psi_{t_{\infty}}^{\eta_{t_{\infty}}}}\l[1-\expp{-\lz\sigma}\r]\cr
\ar=\ar (\psi_{t_{\infty}}^{\eta_{t_{\infty}}})^{-1}(\lz)\cr
\ar=\ar \psi_{t_{\infty}}^{-1}(\lz)-\eta_{t_{\infty}}.
\eeqlb
\qed

 Recall that $\cT_t=\cT\cap [t,\infty)$ for $t\in\cT$.
   According to {\it Step 2} in the proof of   Theorem \ref{MainA},
 for $q\in\cT$, $\Psi_q^{\eta_q}=\{\psi_t^{\eta_q}, t\in\cT_q\}$
 is also an admissible family satisfying (H1-3),
 where $\psi_t^{\eta_q}(\lambda)=\psi_t(\lambda+\eta_q)-\psi_t(\eta_q)$.
 Set $\cT^q_0=\{\theta\geq0, \theta+q\in\cT_q\}$ and $\Psi^q=\{\psi_{\theta+q}^{\eta_q}, \theta\in\cT_0^q\}$.
 Then we have the following corollary.
\begin{corollary}\label{coroTA+}
 For $q\in (t_{\infty}, 0)$, for any nonnegative measurable functional $F$,
\[
\bN^{\Psi}[F(\T_{A+t}, t\in\cT^q_0)|A=q]=\psi'_{q}(\eta_{q})\bN^{\Psi^q}[F(\T_t, t\in\cT^q_0)\sigma_0].
\]
\end{corollary}
\proof Using (\ref{Gir3}) and (\ref{TA}), we have for any nonnegative measurable functional $F$ on $\mT$,
\beqlb\label{taq-etaq-sigma}
\bN^{\Psi}[F(\T_A)|A=q]=\psi'_{q}(\eta_{q})
\mN^{\psi_{q}^{\eta_q}}[F(\T)\sigma].
\eeqlb
Note that $\cT_0^q=\cT_q-q$. Then the desired result follows from the fact that $\{\psi_t, t\in\cT_q\}$
and $\{\psi_t^{\eta_q}, t\in\cT_q\}$ induce the same pruning parameters. \qed

An application of Corollary \ref{coroTA+} is to study the distribution of exit times. Define
\[
A_h=\sup\{t\in\cT; H_{max}(\T_t)>h\},\quad h>0
\]
with the convention $\sup\emptyset=t_{\infty}.$ Then $A_h$ is the last time that the height of the
trees is larger than $h$. The next result is a generalization of Proposition 4.3 in \cite{ADH0a} which
computes the conditional distribution of $A_h$, given $A$.
\begin{proposition}\label{Ah}
 For $t_{\infty}<q<q_0<0$, we have
\beqnn\bN^{\Psi}[A_h>q_0|A=q]
\ar=\ar \frac{\psi'_{q_0}(\eta_q)}{\psi'_{q}(\eta_q)}-\psi'_{q_0}(\eta_q)\psi_q^{\eta_q}(v^{\psi_q^{\eta_q}}(h))
\int_{v^{\psi_q^{\eta_q}}(h)}^{\infty}\frac{dr}{\psi_q^{\eta_q}(r)^2};\cr
\bN^{\Psi}[A_h=A|A=q]
\ar=\ar\psi'_{q}(\eta_q)\psi_q^{\eta_q}(v^{\psi_q^{\eta_q}}(h))
\int_{v^{\psi_q^{\eta_q}}(h)}^{\infty}\frac{dr}{\psi_q^{\eta_q}(r)^2}.
\eeqnn
\end{proposition}
\proof The second equality follows from the fact $\bN^{\Psi}[A_h\geq q|A=q]=1$ and the first equality as $q_0\rar q$.
 We only need to prove the first one. Recall (\ref{Za}).
 We use $\cz_a(\T)$ here  to stress the dependence on $\T$.  Note that
 \beqnn
\bN^{\Psi}[A_h>q_0|A=q]=\bN^{\Psi}[\cz_h(\T_{q_0})>0|A=q]=\bN^{\Psi}[\cz_h({\T_{A+q_0-q}})>0|A=q],
\eeqnn
which, by Corollary \ref{coroTA+}, is equal to
$
\psi_q'(\eta_q)\bN^{\Psi^q}[\ind_{\{\cz_h(\T_{q_0-q})>0\}}\sigma_0].
$
Since for every $t\in\cT_q$, $\psi_t^{\eta_q}$ is subcritical, by (\ref{lemsigma01}), we have
\beqnn
\bN^{\Psi}[A_h>q_0|A=q]\ar=\ar
\psi_q'(\eta_q)\bN^{\Psi^q}\l[\ind_{\{\cz_h(\T_{q_0-q})>0\}}\bN^{\Psi^q}[\sigma_0|\T_{q_0-q}]\r]\cr
\ar=\ar\psi'_{q_0}(\eta_q)\bN^{\Psi^q}\l[\ind_{\{\cz_h(\T_{q_0-q})>0\}}\sigma_{q_0-q}\r]\cr
\ar=\ar\psi'_{q_0}(\eta_q)\mN^{\psi_q^{\eta_q}}\l[\ind_{\{\cz_h>0\}}\sigma\r]\cr
\ar=\ar\psi'_{q_0}(\eta_q)\mN^{\psi_q^{\eta_q}}[\sigma]
-\psi'_{q_0}(\eta_q)\mN^{\psi_q^{\eta_q}}\l[\ind_{\{\cz_h=0\}}\int_0^h\cz_ada\r]\cr
\ar=\ar \frac{\psi'_{q_0}(\eta_q)}{\psi'_{q}(\eta_q)}
-\psi'_{q_0}(\eta_q)\int_0^hda\lim_{\lz\rar0}\mN^{\psi_q^{\eta_q}}\l[\cz_ae^{-\lz\cz_h}\r],
\eeqnn
where we used (\ref{taq-etaq-sigma}) in the last equality.

Recall (\ref{va}), (\ref{ua}) and (\ref{uva}). Then
by (\ref{sigma}) and branching property (iv), conditioning on $\cz_a$, we yield
\[
\lim_{\lz\rar0}\mN^{\psi_q^{\eta_q}}\l[\cz_ae^{-\lz\cz_h}\r]
=\lim_{\lz\rar0}\mN^{\psi_q^{\eta_q}}\l[\cz_ae^{-\cz_au^{\psi_q^{\eta_q}}(h-a,\lz)}\r]
=\frac{\partial}{\partial \lz}u^{\psi_q^{\eta_q}}(a, v^{\psi_q^{\eta_q}}(h-a)).
\]
Since
\begin{align*}
\frac{\partial}{\partial \lz}u^{\psi_q^{\eta_q}}(a, v^{\psi_q^{\eta_q}}(h-a))
&=\frac{\psi_q^{\eta_q}(u^{\psi_q^{\eta_q}}(a, v^{\psi_q^{\eta_q}}(h-a)))}{\psi_q^{\eta_q}(v^{\psi_q^{\eta_q}}(h-a))}\\
&=\frac{\psi_q^{\eta_q}(v^{\psi_q^{\eta_q}}(h))}{\psi_q^{\eta_q}(v^{\psi_q^{\eta_q}}(h-a))^2}
\frac{\partial}{\partial a}v^{\psi_q^{\eta_q}}(h-a).
\end{align*}
Then we conclude that
\beqnn\bN^{\Psi}[A_h>q_0|A=q]\ar=\ar
\frac{\psi'_{q_0}(\eta_q)}{\psi'_{q}(\eta_q)}-\psi'_{q_0}(\eta_q)\int_0^h\frac{\partial}{\partial \lz}
u^{\psi_q^{\eta_q}}(a, v^{\psi_q^{\eta_q}}(h-a))da\cr
\ar=\ar \frac{\psi'_{q_0}(\eta_q)}{\psi'_{q}(\eta_q)}-\psi'_{q_0}(\eta_q)\psi_q^{\eta_q}(v^{\psi_q^{\eta_q}}(h))
\int_{v^{\psi_q^{\eta_q}}(h)}^{\infty}\frac{dr}{\psi_q^{\eta_q}(r)^2}.
\eeqnn
\qed
\begin{remark}
It is easy to see that $\bN^{\Psi}[A_h\geq q]=v^{\psi_q}(h)$.
\end{remark}
\begin{remark}
Recall Theorem \ref{MainA} and Remark \ref{remincreasing},
 an explicit construction of an increasing tree-valued process as that of \cite{ADH0a} may be given which has the
 same distribution as $\{\T_q, q\in \cT_t\}$ under $\bN^{\Psi}$.
 Then by  similar arguments as in \cite{ADH0a} (Theorem 4.6 there),
 one can derive the joint distribution of $(\T_{A_h-}, \T_{A_h})$ (and hence $(\T_{A-}, \T_{A})$).
We left these to the interested readers.
\end{remark}

\section{Tree at the ascension time}\label{sec:ascension}
In this section, we study the representation of the tree at the ascension time. Recall that we shall always assume that
\[
0\in\cT,\quad t_{\infty}<0, \ \mbox{and}\  \psi_0 \ \mbox{is\ critical}.
\]
We first consider an infinite CRT and its pruning. An  infinite  CRT was
 constructed  in  \cite{AD12}  which is the L\'evy CRT conditioned to have infinite height. Before recalling
its construction, we stress that under $\mP_r^\psi$, the root $\emptyset$
belongs  to  $\Br_\infty  $  and  has  mass  $\Delta_0=r$.  We
identify the half real line $[0,+\infty )$ with a real tree denoted by
$\llbracket 0,\infty \llbracket$ with the null mass measure. We denote
by $dx$ the length measure on $\llbracket 0,\infty \llbracket$.
 Let
$\sum_{i\in I_1^*} \dz_{(x^{*}_i, \T^{*,i})}$ and $\sum_{i\in I_2^*} \dz_{(x^{*}_i, \T^{*,i})}$
be two independent Poisson random measures on
$\llbracket 0,\infty \llbracket \times    {\mT}$    with     intensities
$$dx \, 2c\mN^{\psi_0}[d\T] \quad\text{and}\quad dx\int_0^{\infty}lm_0(dl)\mP_l^{\psi_0}(d\T),$$
respectively. The infinite CRT from \cite{AD12}  is defined as
\begin{equation}
   \label{eq:T*}
\ct^*=\lb\emptyset,\infty\lb\circledast_{i\in I_1^*\cup  I_2^*}(x^{*,i}, T^{*,i}).
\end{equation}
We denote by $\mP^{*,\psi_0}(d\ct^*)$ the distribution of $\ct^*$ and $\mE^{*,\psi_0}(d\ct^*)$
the corresponding expectation.
Similarly  to the  setting  in  Section
\ref{sec:prune},   we    consider   on   $\ct^*$   the mark   processes
$M^{\ct^*}_{ske}(dq, dy)$ and $M^{\ct^*}_{node}(dq, dy)$ which are  Poisson random measures
on $\cT_0\times \T^*$ with intensities
\[
\ind_{\{q\in\cT_0\}}\beta_q dq \ell^{\ct^*}(dy)\quad\text{and}\quad\ind_{\{q\in\cT_0\}}\sum_{i\in I_1^*\cup I_2^*}
\sum_{x\in \Br_{\infty}(T^{*,i  })}m_{\Delta_x}(0, dq)\delta_x(dy),
\]
respectively. We identify $x^{*, i}$ as the root of $\ct^{*,i}$. In
particular nodes in $\lb\emptyset,\infty\lb$ with infinite degree will
be charged by $M^{\ct^*}_{node}$. Then set
\[
M^{\ct^*}(dq, dy)=M^{\ct^*}_{ske}(dq, dy)+M^{\ct^*}_{node}(dq, dy).
\]
For every $q\in\cT_0$,  the  pruned tree at  time $q$ is defined as
\[
\ct_q^*=\{x\in\T^*;\ M^{\T^*}([0,q]\times \lb\emptyset,x\lb)=0\},
\]
 with the induced metric, root  $\emptyset$ and mass
measure restricted to $\ct_q^*$.
Our main result in this section is the following theorem and
the proof will be postponed to the end of this section.
\begin{theo}
\label{therep}
  Given $q\in (t_{\infty}, 0)$, if there exists $\bar{q}\in\cT_0$
 such that $\psi_{\bar{q}}(\lz)=\psi_q(\eta_q+\lz)$,
 then conditioned on $\{A=q\}$, $\T_A$ is distributed as $\T_{\bar{q}}^*.$
\end{theo}

\medskip
We remark here that this result is a generalization of the particular setting in \cite{AD12}
or of the discrete case considered in \cite{AP98, ADH12a}. Now
we give some applications of it. Recall
$\cT^q_0=\{\theta\geq0, \theta+q\in\cT_q\}$.
A similar reasoning as Corollary \ref{coroTA+} yields the following corollary.
\begin{corollary}\label{correpb}
Given $q\in (t_{\infty}, 0)$,  if there exists $\bar{q}\in\cT_0$
 such that $\psi_{\bar{q}+t}(\lz)=\psi_{q+t}^{\eta_q}(\lz)$ for all $t\in\cT^{\bar{q}}_0$,
 then conditioned on $\{A=q\}$, $\{\T_{A+t}, t\in\cT^{\bar{q}}_0\}$
 is distributed as $\{\T_{\bar{q}+t}^*,t\in\cT^{\bar{q}}_0\}.$
\end{corollary}
\begin{lemma}
\label{lemrep}Assume that $t_{\infty}\notin \cT_{\infty}$ and for every $t\in(t_{\infty}, 0)$,
there exists $\bar{t}\in \cT$ such that $\psi_{\bar{t}}(\lz)=\psi_{t}(\eta_t+\lz).$
Then $t\rar\bar{t}$ is differentiable and
\[
\frac{d\bar{t}}{dt}
=\frac{\zeta'_t(\eta_t)\psi'_t(\eta_t)-\psi''_t(\eta_t)\zeta_t(\eta_t)}{\zeta'_{\bar{t}}(0)\psi'_t(\eta_t)}
=:\frac{-\gamma_t}{\zeta'_{\bar{t}}(0)},
\quad t\in(t_{\infty}, 0).
\]
\end{lemma}
\proof  It is obvious that  $t\rar\bar{t}$ is differentiable.
By $\psi_{\bar{t}}(\lz)=\psi_t(\eta_t+\lz)$ and (\ref{detat}), we have for all $\lz>0$,
\beqnn
\frac{d\bar{t}}{dt}=\frac{\zeta_t(\eta_t+\lz)\psi'_t(\eta_t)
-\psi'_t(\eta_t+\lz)\zeta_t(\eta_t)}{\zeta_{\bar{t}}(\lz)\psi'_t(\eta_t)}.
\eeqnn
The result follows by taking $\lambda\rar 0$.\qed

Define $\bar{t}_{\infty}=\sup\{\bar{t}; t\in\cT, b_t<0\}$. For $t\in (0, \bar{t}_{\infty})$,
let $\hat{t}$ be the unique negative number such that
$\bar{\hat{t}}=t.$
Let $U$ be a positive ``random'' variable with nonnegative ``density'' with respect to the Lebesgue measure given by
\[
\ind_{\{t\in(0, \bar{t}_{\infty})\}}
\frac{\zeta_{\hat{t}}(\eta_{\hat{t}})\zeta'_t(0)}{\psi'_{{\hat{t}}}(\eta_{\hat{t}})\gamma_{\hat{t}}}.
\]
Assume that $U$ is independent of $\T^*$.
\begin{corollary}
\label{correp}
Suppose that all assumptions in Lemma \ref{lemrep} hold. Then $\T_A$ is distributed under $\bN^{\Psi}$ as $\T_U^*$.
\end{corollary}
\begin{remark}
If $U$ has the same distribution as $A$, then we have $\T_A$ is distributed  as $\T_{\bar{U}}^*.$
\end{remark}
\proof
Recall (\ref{detat}). By Lemma \ref{lemA}, we have the law of $A$ under $\bN^{\Psi}$
has a density with respect to the Lebesgue measure on $\mbb R$ give by
\[
\ind_{\{t\in(t_{\infty}, 0)\}}\frac{\zeta_{t}(\eta_{t})}{\psi'_{t}(\eta_{t})}.
\]
Thus for any nonnegative measurable function $F$ on $\mT$, we deduce from Theorem \ref{therep} that
\beqnn
\bN^{\Psi}[F(\T_A)]\ar=\ar \int_{(t_{\infty}, 0)}\mE^{*, \psi_0}[F(\T^*_{\bar{t}})]
\frac{\zeta_{t}(\eta_{t})}{\psi'_{t}(\eta_{t})}dt\cr
\ar=\ar \int_{(t_{\infty}, 0)}\mE^{*, \psi_0}[F(\T^*_{t})]
\frac{\zeta_{\hat{t}}(\eta_{\hat{t}})}{\psi'_{\hat{t}}(\eta_{\hat{t}})}d\hat{t}\cr
\ar=\ar \int_{(t_{\infty}, 0)}\mE^{*, \psi_0}[F(\T^*_{t})]
\frac{\zeta_{\hat{t}}(\eta_{\hat{t}})}{\psi'_{\hat{t}}(\eta_{\hat{t}})}\frac{d\hat{t}}{d\bar{\hat{t}}}dt
\cr
\ar=\ar \int_{(0, \bar{t}_{\infty})}\mE^{*, \psi_0}[F(\T^*_{t})]
\frac{\zeta_{\hat{t}}(\eta_{\hat{t}})}{\psi'_{\hat{t}}(\eta_{\hat{t}})}
\frac{\zeta'_t(0)}{\gamma_{\hat{t}}}dt\cr
\ar=\ar  \mE^{*, \psi_0}[F(\T_U^*)],
\eeqnn
where the fourth equality follows from Lemma \ref{lemrep}. We have completed the proof.
\qed

 By Corollaries \ref{correpb} and \ref{correp},
 we derive the following result which is a generalization of Corollary 8.2 in \cite{AD12}.
\begin{corollary}\label{correpa}
Suppose $t_{\infty}\notin\cT_{\infty}$ and $[0,\infty)\subset\cT$.
If for every $q\in (t_{\infty}, 0)$, there exists $\bar{q}\in \cT_0$ such that
$\psi_{\bar{q}+t}(\lz)=\psi_{q+t}^{\eta_q}(\lz)$ for all $t\in\cT_0$, then
 $\{\T_{A+t}, t\geq0\}$ is distributed under $\bN^{\Psi}$ as $\{\T_{U+t}^*,t\geq0\}$.
\end{corollary}

In the following we give some examples.
\begin{example}
Recall $\{\psi^{q}(\lz), q\in \Theta^{\psi}\}$ in Example \ref{exaAD12}.
It was assumed in \cite{AD12} that $\psi$ is critical, $\theta_{\infty}:=\inf\Theta^{\psi}\notin\Theta^{\psi}$
and $\theta_{\infty}<0$. So $\{\psi^q, q\in\Theta^{\psi}\}$ satisfies assumptions in Corollary
\ref{correpa}. Then for $t\in(\theta_{\infty}, 0)$, $\eta_t=\bar{t}-t$
and
\[
\zeta_t(\lz)=2c\lz+\int_0^{\infty}(1- \expp{-z\lz})\expp{-zt}zm(dz)=\psi'(t+\lz)-\psi'(t).
\]
Using $\eta_{\hat{t}}=t-\hat{t}$, it is easy to see
\[
\ind_{\{t\in(0, \bar{\theta}_{\infty})\}}
\frac{\zeta_{\hat{t}}(\eta_{\hat{t}})\zeta'_t(0)}{\psi'_{{\hat{t}}}(\eta_{\hat{t}})\gamma_{\hat{t}}}
=\ind_{\{t\in(0, \bar{\theta}_{\infty})\}}\l(1-\frac{\psi'(t)}{\psi'(\hat{t})}\r).
\]
Then we go back to Corollary 8.2 in \cite{AD12}.
\end{example}
\begin{example}
Let $b, c>0$ be two constants. Define $\psi_q({\lz})=qb\lz+c\lz^2$ with $q\in{\mbb R}$ and $\lz\geq0$.
Then $\{\psi_q, q\in\mR\}$ satisfies assumptions in Corollary \ref{correpa}.
In particular, if $b=2c$, we have $\psi_q(\lz)=\psi_0(q+\lz)-\psi_0(q).$
\end{example}
\begin{example}
Recall $\{\psi_q, q\in \cT_-\cup(-\cT_-)\}$ considered in Example \ref{exanobody}.
It is easy to verify that $\{\psi_q, q\in \cT_-\cup(-\cT_-)\}$ satisfies assumptions in Corollary \ref{correp}.
\end{example}

\medskip

\noindent The end of this section is devoted to the proof of Theorem \ref{therep}.

\begin{proposition}\label{proprep}
For any nonnegative measurable functional $F$ on $\mT$ and for every $q\in\cT_0$,
\beqlb\label{therep2}
\psi'_{q}(0)\bN^{\Psi}\l[\sigma_{q}F(\T_{q})\r]
={\mE}^{*,\psi_0}\l[F(\T_{q}^*)\r].
\eeqlb
\end{proposition}
\proof First, recall the Bismut decomposition of L\'evy tree $\T$ along a spine $\lb\emptyset,x\rb$ for $x\in \mbox{Lf}(\T)$. Take the spine as a
subtree, and consider the connected component $(\T^i,i\in I_x)$ of $\T\setminus \lb\emptyset,x\rb$.
Let the branching point $(x_i,i\in I_x)$ be the root. Then $\T=\lb\emptyset,x\rb \circledast_{i\in I_x} (\T^i,x_i)$.
We deduce from Theorem 2.18 in \cite{ADH0a} (or Theorem 7.1 in \cite{AD12}, which were originally proposed in \cite{DL05}) and $\sigma=m^{\T}(\T)$ that,

\beqlb\label{proprep1}\psi'_{q}(0)\bN^{\Psi}\l[\sigma_{q}F(\T_{q})\r]
\ar=\ar\psi'_q(0)\mN^{\psi_q}\l[\sigma F(\T)\r]\cr
\ar=\ar\psi'_q(0)\mN^{\psi_q}\l[\int m^{\T}(dx )F(\lb\emptyset,x\rb \circledast_{i\in I_x} (\T^i,x_i))\r]\cr
\ar=\ar\psi'_q(0)\int_0^{\infty}da
\expp{-\psi'_q(0)a}\mE\l[F\l(\lb\emptyset,a\rb\circledast_{i\in \tilde{I},z_i\leq a}  \tilde{\T}^i\r)\r],
\eeqlb
where under $\mE$, $\sum_{i\in \tilde{I}}\dz_{(z_i, \tilde{\T}^i)}(dz, d\T)$
is a Poisson random measure on $[0, \infty)\times \mT$ with intensity
\[
dz\l(2c\mN^{\psi_q}[d\T]+\int_0^{\infty}lm_q(dl)\mP_l^{\psi_q}(d\T)\r).
\]
For $i\in I_1^*\cup I_2^*$, define
\[
\T_q^{*,i}=\{x\in \T^{*,i}: M^{\T^*}([0,q]\times\lb\emptyset, x\lb)=0\},\quad q\in\cT_0.
\]
With abuse of notation, we have
\beqlb\label{Tqstar}
\T_q^*=\lb\emptyset,\xi\rb\circledast_{i\in I_1^*\cup I_2^*,\; x^{*,i}< \xi}(x^{*,i}, \T_q^{*,i}),
\eeqlb
where
\beqnn
\xi\ar:=\ar\sup\{x\in\lb\emptyset, +\infty\lb: M^{\T^*}([0,q]\times\lb\emptyset, x\lb)=0\}\cr
\ar=\ar\sup\{x\in\lb\emptyset, +\infty\lb: M^{\T^*}_{ske}([0,q]\times\lb\emptyset, x\lb)=0\}\wedge
\inf\{x^{*,i}:  M^{\T^*}_{node}([0,q]\times\{x^{*,i}\})>0\}\cr
\ar=:\ar \xi_1\wedge\xi_2.
\eeqnn

Thanks to (\ref{proprep1}) and (\ref{Tqstar}), it suffices to show that
$\xi$ is exponentially distributed with parameter $\psi'_q(0)$.
Indeed, it is obvious that $\xi_1$ is exponentially distributed with parameter $\int_0^q\beta_{\theta}d\theta$.
By Corollary \ref{coromain1} and the property of Poisson random measure, we have
\[
\sum_{i\in I_2^*}\ind_{\{M^{\T^*}_{node}([0,q]\times\{x^{*,i}\})>0\}}\dz_{(x^{*,i}, \T_q^{*, i})}(dx, d\T)
\]
is a Poisson random measure with intensity
$dx \int_{0}^qd\theta\int_0^{\infty}zn_{\theta}(dz)\mP_z^{\psi_q}(d\T).$
So we  deduce that
$\xi_2$ is exponentially distributed with parameter
$\int_{0}^qd\theta\int_0^{\infty}zn_{\theta}(dz)$.
Hence  $\xi$ is exponentially distributed with parameter
\[
\int_0^q\beta_{\theta}d\theta+ \int_{0}^qd\theta\int_0^{\infty}zn_{\theta}(dz),
\]
 which, by (\ref{b_q}), is just $b_q=\psi'_q(0).$ The result follows.
 \qed

 Now we are in position to prove Theorem \ref{therep}.

\noindent{\bf Proof of Theorem \ref{therep}}:
For any nonnegative measurable function $F$ on $\mT$, by (\ref{TA1}), we have for $q<0$,
\beqlb\label{therep1}
\bN^{\Psi}[F(\T_A)|A=q]\ar=\ar\psi'_{q}(\eta_{q})
\mN^{\psi_{q}}\l[F(\T)\sigma\ind_{\{\sigma<\infty\}}\r]\cr
\ar=\ar \psi'_{q}(\eta_{q})
\mN^{\psi_{q}^{\eta_q}}[F(\T)\sigma],
\eeqlb
where  the last equality follows from (\ref{Gir3}).
Since $\psi_{\bar{q}}(\lz)=\psi_{q}(\eta_{q}+\lz)=\psi_q^{\eta_q}(\lz)$
and $\psi'_q(\eta_q)=\psi'_{\bar{q}}(0)$, an application of Proposition \ref{proprep} yields
\[
\bN^{\Psi}[F(\T_A)|A=q]=\psi'_{\bar{q}}(0)\bN^{\Psi}\l[\sigma_{\bar{q}}F(\T_{\bar{q}})\r]
=\mE^{*,\psi_0}[F(\T_{\bar{q}}^*)].
\]
We have completed the proof. \qed

\bigskip




\end{document}